\documentclass[11pt]{article}
\usepackage{amsthm,amssymb,amsmath, color}
\usepackage{amsfonts}
\usepackage{CJK}
\usepackage{mathrsfs}
\usepackage{amsfonts,amssymb}
\usepackage{amsmath,amscd,mathptm}
\usepackage[all]{xy}
\usepackage{pstricks}
\usepackage{fancyhdr}
\usepackage{wasysym, stmaryrd}
\usepackage[dvipdfm,
            colorlinks=ture,
            citecolor=blue,
            linkcolor=black]{hyperref}

\begin{document}
\sloppy
\newcommand{\dickebox}{{\vrule height5pt width5pt depth0pt}}
\newtheorem{Def}{Definition}[section]
\newtheorem{Bsp}{Example}[section]
\newtheorem{Prop}[Def]{Proposition}
\newtheorem{Theo}[Def]{Theorem}
\newtheorem{Rem}[Def]{Remark}
\newtheorem{Lem}[Def]{Lemma}
\newtheorem{Koro}[Def]{Corollary}
\newcommand{\lra}{\longrightarrow}
\newcommand{\ra}{\rightarrow}
\newcommand{\F}{\mathcal {F}}
\newcommand{\Hom}{{\rm Hom}}
\newcommand{\End}{{\rm End}}
\newcommand{\Ext}{{\rm Ext}}
\newcommand{\Tor}{{\rm Tor}}
\newcommand{\pd}{{\rm proj.dim}}
\newcommand{\inj}{{\rm inj}}
\newcommand{\lgd}{{l.{\rm gl.dim}}}
\newcommand{\gld}{{\rm gl.dim}}
\newcommand{\fd}{{\rm fin.dim}}
\newcommand{\Fd}{{\rm Fin.dim}}
\newcommand{\lfd}{l.{\rm Fin.dim}}
\newcommand{\rfd}{r.{\rm Fin.dim}}
\newcommand{\Mod}{{\rm Mod}}
\newcommand{\Proj}{{\rm Proj}}
\newcommand{\modcat}[1]{#1\mbox{{\rm -mod}}}
\newcommand{\pmodcat}[1]{#1\mbox{{\rm -proj}}}
\newcommand{\Pmodcat}[1]{#1\mbox{{\rm -Proj}}}
\newcommand{\injmodcat}[1]{#1\mbox{{\rm -inj}}}
\newcommand{\E}{{\rm E}_{\mathcal {F}}^{{\rm F},\Phi}}
\newcommand{\X}{ \mathscr{X}_{\mathcal {F}}^{{\rm F},\Phi}}
\newcommand{\Y}{\rm \mathscr{Y}_{\mathcal {F}}^{{\rm F},\Phi}}
\newcommand{\A}{\mathcal {A}}
\newcommand{\C}{\rm \mathscr{C}}
\newcommand{\K}{\rm \mathscr{K}}
\newcommand{\D}{\rm \mathscr{D}}
\newcommand{\opp}{^{\rm op}}
\newcommand{\otimesL}{\otimes^{\rm\bf L}}
\newcommand{\otimesP}{\otimes^{\bullet}}
\newcommand{\rHom}{{\rm\bf R}{\rm Hom}}
\newcommand{\projdim}{\pd}
\newcommand{\stmodcat}[1]{#1\mbox{{\rm -{\underline{mod}}}}}
\newcommand{\Modcat}[1]{#1\mbox{{\rm -Mod}}}
\newcommand{\modcatr}[1]{\mbox{{\rm mod}#1}}
\newcommand{\Modcatr}[1]{\mbox{{\rm Mod}#1}}
\newcommand{\Pmodcatr}[1]{\mbox{{\rm Proj}#1}}

\newcommand{\procat}[1]{#1\mbox{{\rm -proj}}}
\newcommand{\Tr}{\rm Tr}
\newcommand{\add}{{\rm add}}
\newcommand{\Imf}{{\rm Im}}
\newcommand{\Ker}{{\rm Ker}}
\newcommand{\EA}{{\rm E^\Phi_\mathcal {A}}}
\newcommand{\pro}{{\rm pro}}
\newcommand{\Coker}{{\rm Coker}}
\newcommand{\id}{{\rm id}}
\newcommand{\M}{\mathcal {M}}
\newcommand{\Mf}{\rm \mathcal {M}^f}
\newcommand{\rad}{{\rm rad}}
\newcommand{\injdim}{{\rm inj.dim}}

{\large \bf
\begin{center} Derived equivalences in $n$-angulated categories\end{center}}
\medskip
\centerline{\bf Yiping Chen }
\medskip \centerline{School of Mathematics and Statistics, Wuhan University, Wuhan 430072, China}
\medskip \centerline{School of Mathematical Sciences, Beijing Normal University, Beijing 100875, China}

\medskip
\centerline{E-mail: ypchen@whu.edu.cn}

\renewcommand{\thefootnote}{\alph{footnote}}
\setcounter{footnote}{-1} 
\renewcommand{\thefootnote}{\alph{footnote}}
\setcounter{footnote}{-1} 

\abstract{
In this paper, we consider $n$-perforated Yoneda algebras for
$n$-angulated categories, and show that, under some conditions,
$n$-angles induce derived equivalences between the quotient algebras
of $n$-perforated Yoneda algebras. This result generalizes some
results of Hu, K\"{o}nig and Xi. And it also establishes a
connection between higher cluster theory and derived equivalences.
Namely, in a cluster tilting subcategory of a triangulated category,
an Auslander-Reiten $n$-angle implies a derived equivalence between
two quotient algebras. This result can be compared with the fact
that an Auslander-Reiten sequence suggests a derived equivalence
between two algebras which was proved by Hu and Xi.  }

\renewcommand{\thefootnote}{\alph{footnote}}

\section{Introduction}

Derived categories and derived equivalences occur widely in a number
of mathematical fields. For example, algebraic geometry \cite{AAB,
BGG, LM}, differential equation \cite{MS, MK}, the representation
theory of algebras \cite{BK, LS}. In modern representation theory of
finite groups, the famous Abelian defect conjecture of Brou\'{e} is
actually to predicate a derived equivalence between two block
algebras. As is known, derived equivalences preserve many
homological properties of algebras such as the number of simple
modules, the finiteness of global dimension and finitistic
dimension, the algebraic K-theory and Hochschild (co)homological
groups (see \cite{MB, DS, K, R1, R2, PX}). In this sense, derived
equivalences provide us a bridge to compare properties of different
algebras, and are helpful for us to understand some properties of
algebras through the other ones. One of the fundamental problems on
the study of derived equivalences of rings is

{\bf How to construct derived equivalences between rings?}

Richard gave a theoretical solution to this problem which is well
known as the Morita theorem for derived categories \cite{R1} (see
also Keller \cite{K1}). The Richard's theorem for derived categories
is that for two rings $A$ and $B$, the derived categories
$D^b(\Modcat{A})$ and $D^b(\Modcat{B})$ are equivalent as
triangulated categories if and only if there exists a special
complex $T^\bullet$ in $D^b(\Modcat{A})$, called `` tilting
complex'', such that $B$ is the endomorphism ring of $T^\bullet$.
However, it is difficult to construct all tilting complexes
explicitly. And there are so many obstacles to determine the
endomorphism ring of a complex. Consequently, it is necessary to
give a systematic way to construct derived equivalences between
rings.

In order to construct derived equivalences, one strategy is to
develop a practical technique which can produce new derived
equivalences from given ones. In \cite{R1, R2}, Rickard used tensor
product and trivial extension to produce derived equivalences. These
results were generalized by Ladkani in the sense of triangular
matrix ring arising from extension of tilting modules \cite{L1} and
componentwise tensor products \cite{L2}. In \cite{HX2}, Hu and Xi
presented a method to construct new derived equivalences between
these $\Phi$-Auslander Yoneda algebras, or their quotient algebras,
from given almost $\nu$-stable derived equivalences.

Another strategy is trying to construct derived equivalences from
certain sequences. Recently, Hu and Xi introduced $\mathcal
{D}$-split sequences and showed that each $\mathcal {D}$-split
sequence gives rise to a derived equivalence via a tilting module
$\cite{HX1}$. Thus, every Auslander-Reiten sequence is a $\mathcal
{D}$-split sequence and induces a derived equivalence via a
BB-tilting module. This beautiful result presents a relation between
Auslander-Reiten theory and derived equivalences. And later, Hu,
K\"{o}nig and Xi generalized the result in the context of
triangulated categories, adding higher extensions and replacing the
shift functor by any other auto-equivalence of triangulated
categories \cite{HKX}. Note that the derived equivalences are
induced by tilting complexes of length 2. Meanwhile, Ladkani
\cite{L} and Dugas \cite{D} discussed $\mathcal {D}$-split sequences
in the version of mutations of algebras and algebraic triangulated
categories, respectively.

In \cite{GKO}, Geiss, Keller and Oppermann introduced $n$-angulated
categories which occur widely in cluster tilting theory and are
closely related to algebraic geometry and string theory. A natural
question is how to construct derived equivalences in $n$-angulated
categories?

In this paper, we give an affirmative answer to this question. We
construct derived equivalences in the context of $n$-angulated
categories and generalize some results of Hu, K\"{o}nig and Xi in
\cite{HKX}. By the result of Geiss, Keller, Oppermann $\cite{GKO}$,
every $(n-2)$-cluster tilting subcategory which is closed under
$\Sigma^{n-2}$ is an $n$-angulated category. Thus, we can construct
derived equivalences which are induced by tilting complex of
arbitrary length. This result generalizes the main result of Hu,
K\"{o}nig and Xi in \cite{HKX}. At the same time, there is a high
dimensional version of the fact that Auslander-Reiten sequences
suggest a derived equivalence between two algebras which was proved
in \cite{HX1}. Namely, in some cluster tilting subcategory, any
Auslander-Reiten $n$-angle implies a derived equivalence between two
quotient algebras.

In order to describe the main result precisely, we fix some
notations first. Let $R$ be a fixed commutative Artin ring, and let
$k$ be a fixed field. Let $\F$ be an $n$-angulated $R$-category with
suspension functor $\Sigma$, and let $X$ be an object in $\F$.
Suppose that $\F$ has split idempotents. Let $\Phi$ be an admissible
subset of $\mathbb{Z}$. Then we can define $n$-perforated Yoneda
algebra $\E(X):=\oplus_{i\in \Phi}\Hom_{\F}(X, F^iX)$. Its
multiplication is defined in a natural way. The left (right)
$(\add(M), F, \Phi)$-approximation is extension of general
approximation in the sense of Auslander and smal{\o}, adding higher
extension. For more details, we refer readers to section 2. The
objects of $\X(M)$ and $\Y(M)$ satisfy some properties of
orthogonal, i.e.,
$$\X(M):=\{X\in \F\mid \Hom_{\F}(X, F^iM)=0 \mbox{ for all }i\in\Phi/\{0\}\}$$
$$\Y(M):=\{Y\in\F\mid \Hom_{\F}(M, F^iY)=0 \mbox{ for all }i\in\Phi/\{0\}\}.$$

The sets $I$ and $J$ are ideals of $\E(X)$ and $\E(Y)$, respectively
(see section 3 for details). The main result in this paper is the
following:

\begin{Theo}\label{Theo}
  Let $\Phi$ be an admissible subset of $\mathbb{Z}$, and let $\F$
  be an $n$-angulated $R$-category with an auto-equivalence $F$. Suppose
  that $X \stackrel{\alpha_1}\ra M_1\stackrel{\alpha_2}\ra M_2\ra
\cdots \ra M_{n-2}\stackrel{\alpha_{n-1}}\ra Y\stackrel{\alpha_n}\ra
\Sigma X $ is an $n$-angle in $\F$ such that $\alpha_1: X\ra M_1$ is
a left $(\add(M), F, \Phi)$-approximation of $X$ and $\alpha_{n-1}:
M_{n-2}\ra Y$ is a right $(\add(M), F, -\Phi)$-approximation of $Y$.
If $X\in \Y(M)$ and $Y\in \X(M)$, then $\E(X\oplus M)/I$ and
$\E(M\oplus Y )/J$ are derived equivalent.
\end{Theo}

This theorem extends the main result of Hu, K\"{o}nig and Xi in
\cite{GKO}. The following corollary establishes a connection between
higher cluster theory and derived equivalences.

\begin{Koro}
Let $\mathcal{T}$ be a Krull- Schmidt triangulated $k$-category with
shift functor $\Sigma_3$, and let $\mathcal {S}$ be an
$(n-2)$-cluster tilting subcategory of $\mathcal {T}$, which is
closed under $\Sigma^{n-2}_3$. Suppose that
$$X_1\stackrel{\alpha_1}\ra X_2\stackrel{\alpha_2}\ra X_3\ra \cdots\ra X_n$$
is an Auslander-Reiten $n$-angle in $\mathcal{S}$ and $X_1,
X_n\notin \oplus_{i=2}^{n-1}X_i$. Then the two rings
$\End_{\mathcal{S}}(\oplus^{n-1}_{i=1}X_i)/I$ and
$\End_{\mathcal{S}}(\oplus^{n}_{i=2}X_i)/J$ are derived equivalent,
where $I, J$ are defined as in Theorem \ref{Theo}.
\end{Koro}

This paper is organized as follows: In section 2, we make a
preparation for the proof of the main result. We fix some notations
and recall some basic definitions. In section 3, we give the proof
of the main result and deduce some consequences of the main result.
In section 4, we display an example to illustrate our main result.

\section{Preliminaries}

In this section, we will recall some basic definitions and facts
which are needed in our proofs.

\subsection{Notations and conventions}

Throughout this paper, $R$ is a fixed commutative Artin ring with
identity, and $k$ is a fix field.

Let $\cal C$ be an additive category. For an object $X$ in
$\mathcal{C}$, we denote by $\add(X)$ the full subcategory of $\cal
C$ consisting of all direct summands of finite direct sums of $X$.
For two morphisms $f:X\rightarrow Y$ and $g:Y\rightarrow Z$ in $\cal
C$, we write $fg$ for their composition which is a morphism from $X$
to $Z$. For two functors $F:\mathcal{C}\rightarrow \mathcal{D}$ and
$G:\mathcal{D}\rightarrow\mathcal{E}$, we write $GF$ for the
composition instead of $FG$.

Let $\mathcal {C}$ be an additive category with an endo-functor $F:
\mathcal {C}\ra \mathcal {C}$. Let $\mathcal {D}$ be a full
subcategory of $\mathcal {C}$, and let $\Phi$ be a non-empty subset
of $\mathbb{N}$. If $F$ has an inverse, then $\Phi$ can be chosen to
be a subset of $\mathbb{Z}$. Let $X$ be a object of $\mathcal {C}$.
A morphism $f: X\ra D$ in $\mathcal {C}$ is called a {\em left
cohomological $\mathcal {D}$-approximation} of $X$ with respect to
$(F,\Phi)$ (or left $(\mathcal {D}, F, \Phi)$-approximation of $X$)
if $D\in \mathcal {D}$, and for any morphism $g: X\ra F^i(D')$ with
$D'\in \mathcal {D}$ and $i\in\Phi$, there is a morphism $g': D\ra
F^i(D')$ such that $g=fg'$. Note that $F^0=id_\mathcal {C}$. Dually,
we have the notion of {\em right cohomological $\mathcal
{D}$-approximation} of $X$ (or right $(\mathcal {D}, F,
\Phi)$-approximation of $X$) if for any $i\in\Phi$ and any morphism
$g: F^iD'\ra X$ with $D'\in\mathcal {D}$, there is a morphism $g':
F^iD'\ra D$ such that $g=g'f$ (see \cite{HKX}). In particular, if
$\Phi=\{0\}$, then left (resp., right)-$(\mathcal{D}, F,
\Phi)$-approximation of $X$ is left (resp., right)
$\mathcal{D}$-approximation of $X$. The subcategory $\mathcal{D}$ is
called {\em contravariantly finite} subcategory of $\mathcal{C}$ if
any object $Y$ in $\mathcal{C}$ has a right
$\mathcal{D}$-approximation. Dually, a covariantly finite
subcategory of $\mathcal{C}$ is defined. The subcategory
$\mathcal{D}$ is called {\em functorially finite} of $\mathcal{C}$
if $\mathcal{D}$ is contravariantly finite and covariantly finite in
$\mathcal{C}$. We denote by $J_{\mathcal{C}}$ the Jacobson radical
of $\mathcal{C}$. Let $f\in \Hom_{\mathcal{C}}(X, Y)$ be a morphism.
We call $f$ a {\em sink map} of $Y$ if $f$ satisfies the following
conditions: $(1)$ if $g:{}X\ra X$ satisfies $gf=f$, then $g$ is an
automorphism. $(2)$ $f\in J_{\mathcal{C}}$ and
$$\Hom_{\mathcal{C}}(-,X)\stackrel{. f}\ra J_{\mathcal{C}}(-, Y)\ra 0$$
is exact as functors on $\mathcal{C}$. Dually, a source map is
defined (see \cite{IY}).

Given an $R$-algebra $A$, we denote the opposite algebra of $A$ by
$A^{op}$. By an $A$-module we mean a unitary left $A$-module; the
category of all (resp., finitely generated) $A$-modules is denoted
by $A$-Mod (resp., $\modcat{A}$), the full subcategory of $A$-Mod
consisting of all (resp., finitely generated) projective modules is
denoted by $A$-Proj (resp., $\pmodcat{A}$). Similarly, the full
subcategory of $A$-Mod consisting of all (resp., finitely generated)
injective $A$-modules is denoted by $A$-Inj (resp., $A$-inj). An
algebra $A$ is called an Artin $R$-algebra if $A$ is finitely
generated as an $R$-module. Let $A$ be an Artin $R$-algebra, we
denote by $D$ the usual duality on $\modcat{A}$. The functor
$\nu_A:=D\Hom_A(-,{}_AA):{}\pmodcat{A}\ra \injmodcat{A}$ is Nakayama
functor. We denote the syzygy functor by $\Omega$. Namely, for an
$A$-module, we denote the first syzygy of $M$ by $\Omega_A(M)$. The
stable category $A$-\underline{mod} is a quotient category of
$\modcat{A}$. The objects of $A$-\underline{mod} are the objects of
$\modcat{A}$. Let $X, Y$ be in $\modcat{A}$. The homomorphism set
$\underline{\Hom}(X, Y)$ is $\Hom(X, Y)$ modulo the submodule
generated by homomorphism which can factorize through some
projective $A$-module.

Let $A$ be an Artin algebra. A {\em complex} $X^\bullet=(X^i,
d_X^i)$ of $A$-modules is a sequence of $A$-modules and $A$-module
homomorphisms $d^i_X: X^i\ra X^{i+1}$ such that $d^i_Xd^{i+1}_X=0$
for all $i\in \mathbb{Z}$. A {\em morphism} $f^\bullet: X^\bullet\ra
Y^\bullet$ between two complexes $X^\bullet$ and $Y^\bullet$ is a
collection of homomorphisms $f^i: X^i\ra Y^i$ of $A$-modules such
that $f^id^i_Y=d_X^if^{i+1}$. The morphism $f^\bullet$ is said to be
{\em null-homotopic} if there exists a homomorphism $h^i: X^i\ra
Y^{i-1}$ such that $f^i=d_X^ih^{i+1}+h^id_Y^{i-1}$ for all
$i\in\mathbb{Z}$. A complex $X^\bullet$ is called {\em bounded
below} if $X^i=0$ for all but finitely many $i<0$, {\em bounded
above} if $X^i=0$ for all but finitely many $i>0$, and {\em bounded}
if $X^\bullet$ is bounded below and above. We denote by $C(A)$
(resp., $C(\Modcat{A})$) the category of complexes of finitely
generated (resp., all) $A$-modules. The homotopy category $K(A)$ is
quotient category of $C(A)$ modulo the ideals generated by
null-homotopic morphisms. We denote the derived category of
$\modcat{A}$ by $D(A)$ which is the quotient category of $K(A)$ with
respect to the subcategory of $K(A)$ consisting of all the acyclic
complexes. The full subcategory of $K(A)$ and $D(A)$ consisting of
bounded complexes over $\modcat{A}$ is denoted by $K^b(A)$ and
$D^b(A)$, respectively. We denoted by $C^+(A)$ the category of
complexes of bounded below, and by $K^+(A)$ the homotopy category of
$C^+(A)$. The full subcategory of $D(A)$ consisting of bounded below
complexes is denoted by $D^+(A)$. Similarly, we have the category
$C^-(A)$ of complexes bounded above, the homotopy category $K^-(A)$
of $C^-(A)$ and the derived category $D^-(A)$ of $C^-(A)$. If we
focus on the category of left $A$-modules, then we have the homotopy
category $K(\Modcat{A})$ of $C(\Modcat{A})$ and the derived category
$D(\Modcat{A})$ of $C(\Modcat{A})$. 
Suppose that $X^\bullet=(X^i, d_X^i)$ and $Y^\bullet=(Y^i, d_Y^i)$
are two complexes. We define the {\em direct sum} of $X^\bullet$ and
$Y^\bullet$ by the complex $Z^\bullet=(Z^i, d_Z^i)$ such that
$Z^i=X^i\oplus Y^i$ and $d_Z^i=\begin{pmatrix}d_{X}^i&0\\0&d_Y^i
\end{pmatrix}:{}X^i\oplus Y^i\ra X^{i+1}\oplus Y^{i+1}$. The complex
$X^\bullet$ and the complex $Y^\bullet$ are called the {\em direct
summands} of $Z^\bullet$.

\medskip

The following result, due to Rickard (see \cite[Theorem 6.4]{R1}),
may be called the Morita theorem of derived categories.

\begin{Lem}$\cite{R1}$\label{R}
Let $\Lambda$ and $\Gamma$ be two rings. The following conditions
are equivalent:

$(1)$ $K^-(\pmodcat{\Lambda})$ and $K^-(\pmodcat{\Gamma})$ are
equivalent as triangulated categories;

$(2)$ $D^b(\Modcat{\Lambda})$ and $D^b(\Modcat{\Gamma})$ are
equivalent as triangulated categories;

$(3)$ $K^b(\Pmodcat{\Lambda})$ and $K^b(\Pmodcat{\Gamma})$ are
equivalent as triangulated categories;

$(4)$ $K^b(\pmodcat{\Lambda})$ and $K^b(\pmodcat{\Gamma})$ are
equivalent as triangulated categories;

$(5)$ $\Gamma$ is isomorphic to $\End(T^\bullet)$, where $T^\bullet$
is a complex in $K^b(\pmodcat{\Lambda})$ satisfying:

\quad $(a)$ $T^\bullet$ is self-orthogonal, that is,
$\Hom_{K^b(\pmodcat{\Lambda})}(T^\bullet,T^\bullet[i])=0$ for all
$i\neq 0$,

\quad $(b)$ add$(T^\bullet)$ generates $K^b(\pmodcat{\Lambda})$ as a
triangulated category.
\end{Lem}

Two rings $\Lambda$ and $\Gamma$ are called \emph{derived
equivalent} if the above conditions (1)-(5) are satisfied. A complex
$T^\bullet$ in $K^b(\pmodcat{\Lambda})$ as above is called a
\emph{tilting complex} over $\Lambda$. It is also equivalent to say
that the two rings $\Lambda$ and $\Gamma$ are derived equivalent if
and only if there exists a complex $X^\bullet$ in
$D(\Modcat{\Lambda})$, isomorphic to a complex in
$K^{b}(\pmodcat{\Lambda})$ which satisfies [Lemma \ref{R}(5), (a)
and (b)], such that the two rings $\Gamma$ and
$\End_{D(\Modcat{\Lambda})}(X^\bullet)$ are isomorphic. In
particular, if the tilting complex $T^\bullet$ is isomorphic to a
module $T$ in $D^b(\Lambda)$, then $T$ is called {\em tilting
module}.



\medskip
\subsection{The $n$-angulated categories}

In this part, we will recall the definition and some properties of
$n$-angulated categories which are proposed by Geiss, Keller and
Oppermann in $\cite{GKO}$. For the convenience of the reader, we
repeat the relevant material from \cite{GKO}.

Suppose that $\F$ is an additive category with an automorphism
$\Sigma$, and $n$ \;$(\geq 3)$ is an integer. A sequence of objects
and morphisms in $\F$ of the form
$$X_\bullet:=X_1\stackrel{\alpha_1}\ra
X_2\stackrel{\alpha_2}\ra \cdots\stackrel{\alpha_{n-1}}\ra
X_n\stackrel{\alpha_n}\ra \Sigma X_1$$ is called an {\em
$n$-$\Sigma$-sequence}. An $n$-$\Sigma$-sequence $X_\bullet$ is
called {\em exact} if the following sequence of $\mathbb{Z}$-modules
$$\Hom_{\F}(Y,X_\bullet):  \cdots\ra\Hom_{\F}(Y,X_1)\ra \Hom_{\F}(Y,X_2)
\ra\cdots\ra\Hom_{\F}(Y,X_n)\ra\cdots$$ is exact for every object
$Y\in \F$. The left rotation of $X_\bullet$ is the following
$n$-$\Sigma$-sequence
$$X_\bullet[1]:=(X_2\stackrel{\alpha_2}\ra X_3\stackrel{\alpha_3}\ra \cdots
\stackrel{\alpha_n}\ra \Sigma
X_1\stackrel{(-1)^n\Sigma\alpha_1}\ra\Sigma X_2).$$ Similarly, the
right rotation of $X_\bullet$ is the $n$-$\Sigma$-sequence
$$X_\bullet[-1]:=(\Sigma^{-1}X_n\stackrel{(-1)^n\Sigma^{-1}\alpha_n}\ra X_1\stackrel{\alpha_1}\ra\cdots\stackrel{\alpha_{n-2}}\ra X_n).$$

An $n$-$\Sigma$-sequence of the form
$(TX)_\bullet:=(X\stackrel{1_X}\ra X\ra 0\ra\cdots\ra 0\ra\Sigma X)$
for $X\in \F$, or its rotation is called {\em trivial}.  A {\em
morphism} of two $n$-$\Sigma$-sequences is given by a sequence of
morphisms $\varphi=(\varphi_1,\varphi_2,\cdots,\varphi_n)$ in $\F$
such that the following diagram commutates:

$$\xymatrix{
  X_1 \ar[d]_{\varphi_1} \ar[r]^{\alpha_1} & X_2 \ar[d]_{\varphi_2} \ar[r]^{\alpha_2}
  & X_3 \ar[d]_{\varphi_3}\ar[r]  & \cdots \ar[r] & X_n \ar[d]_{\varphi_n} \ar[r]^{\alpha_n} & \Sigma X_1 \ar[d]^{\Sigma\varphi_1} \\
  Y_1 \ar[r]^{\beta_1} & Y_2 \ar[r]^{\beta_2} & Y_3\ar[r]  & \cdots\ar[r]  & Y_n\ar[r]^{\beta_n} & \Sigma Y_n.   }$$

The morphism $\varphi$ is called a {\em weak isomorphism} if
$\varphi_i$ and $\varphi_{i+1}$ are isomorphisms, where $1\leq i\leq
n$, and $\varphi_{n+1}$ is denoted by $\Sigma\varphi_1$. Two
$n$-$\Sigma$-sequences $X^1_\bullet$ and $X^n_\bullet$ are called
{\em weakly isomorphic} if there is a chain of
$n$-$\Sigma$-sequences
$$X^1_\bullet-X^2_\bullet-\cdots-X^{n-1}_\bullet-X^n_\bullet$$
satisfying that there is a weak isomorphism between $X^i_\bullet$
and $X^{i+1}_\bullet$ for $1\leq i\leq n-1$.

\begin{Def}$(\cite{GKO})$
\label{def} A collection $\pentagon$ of $n$-$\Sigma$-sequences is
called a ({\em pre}-) $n$-angulation of $(\F,\Sigma)$ and its
elements {\em n-angles} if $\pentagon$ fulfills the following
conditions:

\begin{enumerate}
\item
\begin{enumerate}
\item $\pentagon$ is closed under direct sums
and under taking summands.
\item For all $X\in \F$, the trivial $n$-$\Sigma$-sequence
$(TX)_\bullet$ belongs to $\pentagon$.
\item For each morphism $\alpha_1: X_1\ra X_2$ is $\F$, there
exists an $n$-angle starting with $\alpha_1$.
\end{enumerate}
\item An $n$-$\Sigma$-sequence $X_\bullet$ belongs to $\pentagon$
if and only if $X_\bullet[1]$ belongs to $\pentagon$.
\item Each commutative diagram

$$\xymatrix{
  X_1 \ar[d]_{\varphi_1} \ar[r]^{\alpha_1} & X_2 \ar[d]_{\varphi_2} \ar[r]^{\alpha_2}
  & X_3 \ar[r]  & \cdots \ar[r] & X_n  \ar[r]^{\alpha_n} & \Sigma X_1 \ar[d]^{\Sigma\varphi_1} \\
  Y_1 \ar[r]^{\beta_1} & Y_2 \ar[r]^{\beta_2} & Y_3\ar[r]  & \cdots\ar[r]  & Y_n\ar[r]^{\beta_n} & \Sigma Y_n   }$$
with rows in $\pentagon$ can be completed to a morphism of
$n$-$\Sigma$-sequences.
\end{enumerate}

Moreover, if $\pentagon$ fulfills the following condition, it is
called an n-angulation of $(\F,\Sigma)$:

4.  In the situation of 3 the morphisms $\varphi_3, \varphi_4,
\cdots, \varphi_n$ can be chosen such that the cone
$C(\varphi_\bullet)$:

$$X_2\oplus Y_1 \stackrel{\begin{pmatrix}
  -\alpha_2&0\\
  \varphi_2&\beta_1
\end{pmatrix}}\lra X_3\oplus Y_2\stackrel{\begin{pmatrix}
  -\alpha_3&0\\
  \varphi_3&\beta_2
\end{pmatrix}}\lra \cdots\stackrel{\begin{pmatrix}
  -\alpha_n&0\\
  \varphi_n&\beta_{n-1}
\end{pmatrix}}\lra \Sigma X_1 \oplus Y_n\stackrel{\begin{pmatrix}
  -\Sigma \alpha_1&0\\
  \Sigma\varphi_1&\beta_n
\end{pmatrix}}\lra \Sigma X_2\oplus \Sigma Y_1$$
belongs to $\pentagon$.
\end{Def}

\medskip

\begin{Def}
Suppose that $(\F,\Sigma,\pentagon)$ and $(\F',\Sigma',\pentagon')$
are two $n$-angulated categories. An additive functor $F: \F\ra\F'$
is called {\em $n$-angle functor} if $F(\pentagon)=\pentagon'$,
i.e., there exists an invertible natural transformation $\xi:
F\Sigma\ra\Sigma'F$ such that
$(FX_1,FX_2,\cdots,FX_n,F\alpha_1,F\alpha_2,\cdots,F\alpha_n\xi_{X_1})$
is in $\pentagon'$ for
$(X_1,X_2,\cdots,X_n,\alpha_1,\alpha_2,\cdots,\alpha_n)$ in
$\pentagon$. Moreover, if $F$ is an equivalence of categories, then
$F$ is called {\em $n$-angle equivalence}.
\end{Def}

\noindent{\bf Remark.} If $n=3$, then $F$ is well-known as triangle
functor.

\bigskip

In $\cite{GKO}$, Geiss, Keller and Oppermann show how to construct
$n$-angulated categories inside triangulated categories.

\begin{Bsp}\label{theoGKO}$\cite{GKO}$ Let $\mathcal {T}$ be a triangulated
category with an $(n-2)$-cluster tilting subcatgory $\F$, which is
closed under $\Sigma^{n-2}_3$, where $\Sigma_3$ denotes the
suspension in $\mathcal {T}$. Then $(\F, \Sigma^{n-2}_3, \pentagon)$
is an $n$-angulated category, where $\pentagon$ is the class of all
sequences

$$X_1\stackrel{\alpha_1}\ra X_2\stackrel{\alpha_2}\ra \cdots\stackrel{\alpha_{n-1}}\ra X_{n}\stackrel{\alpha_n}\ra \Sigma^{n-2}_3 X_1$$
such that there exists a diagram

$$\xymatrix@R=5mm{
   & X_2 \ar[dr]\ar[rr]^{\alpha_2} && X_3 \ar[dr] \ar[rr] && X_4 \ar[dr]\ar[r]& \cdots&\cdots\ar[r]&X_{n-1} \ar[dr]^{\alpha_{n-1}}& \\
 X_1\ar[ur]^{\alpha_1}&&X_{2.5}\ar[ur] \ar@{-->}[ll]&&X_{3.5}\ar[ur] \ar@{-->}[ll]&&\cdots\ar@{-->}[ll]&X_{n-1.5}\ar@{-->}[l] \ar[ur]&&X_n\ar@{-->}[ll]}$$
with $X_i\in \mathcal {T}$ for $i\notin \mathbb{Z}$, such that all
oriented triangles are triangles in $\mathcal {T}$, all non-oriented
triangles commute, and $\alpha_n$ is the composition along the lower
edge of the diagram.
\end{Bsp}

In order to prove the main result, we should prove the following
lemma.

\begin{Lem}\label{lem1.7}
Let $(\F,\Sigma,\pentagon)$ be a pre-$n$-angulated category.

For $2\leq i< n$. Each commutative diagram

$\xymatrix{
  X_1 \ar[d]_{\varphi_1} \ar[r]^{\alpha_1} & X_2 \ar[d]_{\varphi_2}\ar[r]^{\alpha_2} & \cdots \ar[r]^{} & X_i \ar[d]_{\varphi_i} \ar[r]^{} & X_{i+1} \ar@{-->}[d]_{\varphi_{i+1}} \ar[r]^{} & \cdots  \ar[r]^{} & X_n \ar[r]^{\alpha_n}\ar@{-->}[d]_{\varphi_n} & \Sigma X_1\ar[d]_{\Sigma \varphi_1}  \\
 Y_1 \ar[r]^{\beta_1} & Y_2 \ar[r]^{\beta_2} & \cdots \ar[r]^{} & Y_i \ar[r]^{} & Y_{i+1} \ar[r]^{} & \cdots \ar[r]^{} & Y_n \ar[r]^{\beta_n} & \Sigma Y_n
 }$\\
with rows in $\pentagon$ can be completed to a morphism of
$n$-$\Sigma$-sequences.
\end{Lem}

{\bf Proof.} The proof is similar with \cite[Lemma
2.3]{GKO}.$\square$

Suppose that $\F$ has split idempotents. If we denote this lemma by
$(3')$, then we can modify the definition of pre-$n$-angulated
category. That is, a collection $\pentagon$ of
$n$-$\Sigma$-sequences is called a pre-$n$-angulation of $(\F,
\Sigma)$ if $\pentagon$ satisfies the following conditions:
$(1(a)-1(c)), 2, 3'$. It is easy to prove that the two cases of
definition are equivalent. However, the change is vital for the
proof of the main result.

\medskip



\medskip
\subsection {Admissible subsets and $n$-perforated Yoneda algebras}
In this part, we will introduce a new class of algebras which are
called $n$-perforated Yoneda algebras.

 Let ${\mathbb N}=\{0,1, 2, \cdots\}$ be the set of natural
numbers, and let $\mathbb Z$ be the set of all integers. For a
natural number $n$ or infinity, let ${\mathbb N}_n :=\{i\in {\mathbb
N}\mid 0\le i< n+1\}$.

Recall from \cite{HX2} that a subset $\Phi$ of $\mathbb{Z}$
containing $0$ is called an {\em admissible subset} of $\mathbb{Z}$
if the following condition is satisfied:

{\it If $i, j$ and $k$ are in $\Phi$ such that $i+j+k\in\Phi$, then
$i+j\in\Phi$ if and only if $j+k\in\Phi$.}

\medskip
Any subset $\{0,i,j\}$ of $\mathbb N$ is an admissible subset of
$\mathbb Z$. Moreover, for any subset $\Phi$ of $\mathbb N$
containing zero and for any positive integer $m\ge 3$, the set
$\{x^m\mid x\in \Phi\}$ is admissible in $\mathbb Z$. The
intersection of a family of admissible subsets of $\mathbb{N}$ is
admissible (for more examples, see \cite{HX2}). Nevertheless, not
every subset of $\mathbb N$ containing zero is admissible. Note that
$\Phi^2$ is not necessary admissible in $\mathbb{N}$ even if $\Phi$
is an admissible subset of $\mathbb{N}$. For instance, $\{0,1,2,4\}$
is not admissible. In fact, this is the `smallest' non-admissible
subset of $\mathbb N$. For more details, we refer reader to
\cite{HX2}.

Admissible sets were used to define the $\Phi$-Auslander Yoneda
algebras in $\cite{HX2}$ and the perforated Yoneda algebra in
$\cite{HKX}$, if we restrict to the case of an object in a
triangulated category. However, in this paper, we will restrict to
the case of objects in an $n$-angulated category.

The following is the most natural generalization of perforated
Yoneda algebra, proposed by Hu, K\"{o}nig and Xi in \cite{HKX}, for
$n$-angulated categories.

Let $\Phi$ be an admissible subset of $\mathbb{Z}$, and let $\F$ be
an $n$-angulated $R$-category with suspension functor $\Sigma$.
Suppose that $F$ is an $n$-angle functor from $\F$ to $\F$. Note
that $F^i=0$ for $i<0$ if the quasi-inverse of $F$ does not exist.
Consider the $(\Phi, F)$-orbit category ${\F}^{F, \Phi}$, the
extension of orbit category, whose object are the objects of $\F$.
Suppose that $X$ and $Y$ are two objects in ${\F}^{F, \Phi}$, the
homomorphism set in ${\F}^{F, \Phi}$ is defined as follows:
$$\Hom_{{\F}^{F, \Phi}}(X, Y):=\bigoplus_{i\in \Phi}\Hom_{\F}(X, F^i Y)\in \Modcat{R}$$
and the composition is defined in an obvious way. Since $\Phi$ is
admissible, the $(\Phi, F)$-orbit category $\mathcal {F}^{F, \Phi}$
is an additive $R$-category. Let $X, Y$ be objects in ${\F}^{F,
\Phi}$. Thus, $\Hom_{{\F}^{F, \Phi}}(X, X)$ is an $R$-algebra. It is
called the {\em n-perforated Yoneda algebra} of $X$ with respect to
$F$, and denoted by $\E(X)$. $\Hom_{{\F}^{F, \Phi}}(X, Y)$ is a
$\E(X)$-$\E(Y)$-bimodule. For convenience, we denote $\Hom_{{\F}^{F,
\Phi}}(X, Y)$ by $\E(X, Y)$.

The following lemma, which was essentially taken form \cite[Lemma
3.5]{HX2}, \cite[Lemma 2.2]{HKX}, describes the basic properties of
the algebra $\E(X)$ where $X$ is an object in the $n$-angulated
$R$-category $\F$, which can also be verified directly.

\begin{Lem}\label{lem1.0}
Let $\F$ be an $n$-angle $R$-category with an $n$-angle endo-functor
$F$, and let $U$ be an object in $\F$. Suppose that $U_1, U_2, U_3$
are in $\add(U)$, and that $\Phi$ is an admissible subset of
$\mathbb{Z}$. Then

$(1)$ There is a natural isomorphism

$$\mu: \E(U_1,U_2)\ra \Hom_{\E(U)}(\E(U,U_1),\E(U,U_2)),$$
which sends $x\in \E(U_1,U_2)$ to the morphism $a\mapsto ax$ for
$a\in \E(U,U_1)$. Moreover, if $x\in \E(U_1,U_2)$ and $y\in
\E(U_2,U_3)$, then $\mu(xy)=\mu(x)\mu(y)$.

$(2)$ The functor $\E(U,-): add(U)\ra \pmodcat{\E(U)}$ is faithful.

$(3)$ If $\Hom_{\F}(U_1,F^iU_2)=0$ for all $i\in
\Phi\setminus\{0\}$, then the functor $\E(U,-)$ induces an
isomorphism of $R$-modules:

$$\E(U,-): \Hom_{\F}(U_1,U_2)\ra \Hom_{\E(U)}(\E(U,U_1),\E(U,U_2)).$$
\end{Lem}


\medskip

\section{Proof of the main result}

In this section, we will construct derived equivalences from an
$n$-angle. Firstly, we will prove Theorem $\ref{Theo}$. Secondly, we
will derive some consequences form the main result.

Let $\F$ be an $n$-angulated category with suspension functor
$\Sigma$, and let $\pentagon$ be an $n$-angulation of $(\F,
\Sigma)$. Suppose that $\F$ has split idempotent and the functor
$F:{}\F\ra\F$ is an $n$-angle functor. Since $F$ is an $n$-angulated
category, there is a natural isomorphism $\delta: F\Sigma\ra \Sigma
F$ associated with $F$. We denote the isomorphism $F^i(\Sigma^j
X)\ra \Sigma^j(F^iX)$ by $\delta(F,i,X,j)$. Note that there is an
inclusion $\iota: \Hom_{\mathcal {F}}(X, Y)\ra \E(X, Y)$. Given a
morphism $f\in\Hom_{\mathcal {F}}(X, Y)$, $\iota(f)$ is an element
of $\E(X,Y)$ concentrated in degree $0$. For convenience, we denote
$\iota(f)$ by $\underline{f}$.

Set
$$X \stackrel{\alpha_1}\ra M_1\stackrel{\alpha_2}\ra M_2\ra \cdots
\ra M_{n-2}\stackrel{\alpha_{n-1}}\ra Y\stackrel{\alpha_n}\ra \Sigma
X $$ be an $n$-angle in $\pentagon$.

For simplicity, we denote $\oplus_{i=1}^{n-2}M_i$ by $M$ and write
$V, W$ instead of $X\oplus M, M\oplus Y$, respectively. Thus, we can
get $M_i\in \add (M)$ for $i=1, 2, \cdots, n-2$.

Since the direct sum of two $n$-angles is still an $n$-angle, there
are two $n$-angles
$$X \stackrel{\alpha_1}\ra M_1\stackrel{\alpha_2}\ra M_2\ra \cdots\ra
M_{n-3}\stackrel{\overline{\alpha_{n-2}}}\ra M_{n-2}\oplus M
\stackrel{\overline{\alpha_{n-1}}}\ra
W\stackrel{\overline{\alpha_n}}\ra \Sigma X$$
$$\Sigma^{-1}Y\stackrel{(-1)^n\Sigma^{-1}\widetilde{\alpha_n}}\lra
V\stackrel{\widetilde{\alpha_1}}\ra M_1\oplus
M\stackrel{\widetilde{\alpha_2}}\ra
\cdots\stackrel{\alpha_{n-2}}\lra M_{n-2}\stackrel{\alpha_{n-1}}
\lra Y $$

We define
$$\overline{\alpha_{n-2}}:= (\alpha_{n-2}, 0): M_{n-3}\ra
M_{n-2}\oplus M\quad\quad \overline{\alpha_{n-1}}:= \begin{pmatrix}
  0&\alpha_{n-1}\\
  1&0
\end{pmatrix}: M_{n-2}\oplus M\ra M\oplus Y$$
$$\overline{\alpha_n}:= \begin{pmatrix}
0\\
\alpha_n
\end{pmatrix}: M\oplus Y\ra \Sigma X\quad\quad \widetilde{\alpha_1}:=
\begin{pmatrix} \alpha_1&0\\
0&1
\end{pmatrix}:X\oplus M\ra M_1\oplus M$$
$$\widetilde{\alpha_2}:= \begin{pmatrix}\alpha_2\\0
\end{pmatrix}: M_1\oplus M\ra M_2 \quad\quad
\widetilde{\alpha_n}:=\begin{pmatrix} \alpha_n&0
\end{pmatrix}: Y\ra \Sigma V$$

For a subset $\Phi$ of $\mathbb{Z}$, we define $-\Phi:= \{-x\mid
x\in \Phi\}$ and
$$\X(M):= \{X\in \F  \mid  \Hom_{\F}(X, F^iM)=0  \mbox{ for all }
 i\in \Phi \setminus \{0\}\},$$
$$\Y(M):=\{Y\in \F \mid \Hom_{\F}(M, F^iY)=0  \mbox{ for all }
 i\in \Phi \setminus \{0\}\}.$$
$$\begin{array}{ll}I:=\{x=(x_i)\in \E(X\oplus M)\mid& x_i=0\mbox{ for }0\neq i \in \Phi,
\\&x_0 \mbox{ factorizes through add(M) and }
\Sigma^{-1}\widetilde{\alpha_n}\},\end{array}$$
$$\begin{array}{ll}J:=\{y=(y_i)\in \E(M\oplus Y)\mid&y_i=0 \mbox{ for }0\neq i\in \Phi,
\\&y_0 \mbox{ factorizes through add(M) and } \overline{\alpha_n}\}.\end{array}$$

\medskip

In order to prove Theorem \ref{Theo}, we prove the following lemmas.

\begin{Lem}\label{lem1.1}
The sets $I$ and $J$ are ideals of $\E(V)$ and $\E(W)$,
respectively.
\end{Lem}

{\bf Proof.} It is easily seen that the set $I$ is closed under
addition. By the definition of $I$, we can write $x_0=uv$ for $u:
V\ra M'$ and $v: M'\ra V$, where $M'$ is an object in add($M$), and
$x_0=s(\Sigma^{-1}\widetilde{\alpha_n})$ for a morphism $s: V\ra
\Sigma^{-1}Y$. Suppose $x=(x_i)_{i\in\Phi}\in I,
y=(y_i)_{i\in\Phi}\in\E(V)$. In order to prove that the set $I$ is
an ideal of $\E(V)$, it suffices to prove that
$xy=(x_0y_i)_{i\in\Phi}\in I, yx=(y_iF^i(x_0))_{i\in \Phi}\in I$.

It is clear that $x_0y_0$ factorizes through
$\Sigma^{-1}\widetilde{\alpha_n}$ and some object in add($M$). Set
$0\neq i\in\Phi$. Note that $\widetilde{\alpha_1}: V\ra M_1\oplus M$
is a left (add$(M)$, F, $\Phi$)-approximation of $V$. Thus, for
given $y_i: V\ra F^iV$, there is a morphism $z_i: M_1\oplus M\ra
F^i(M_1\oplus M)$ such that
$\widetilde{\alpha_1}z_i=y_iF^i(\widetilde{\alpha_1})$. Since $F$ is
an $n$-angle functor, there is a commutative diagram between two
$n$-angles.
$$\xymatrix{
  \Sigma^{-1}Y \ar[d]_{} \ar[r]^{(-1)^n\Sigma^{-1}\widetilde{\alpha_n}}
  & V \ar[d]_{y_i} \ar[r]^{\widetilde{\alpha_1}} & M_1\oplus M \ar[d]_{z_i}
  \ar[r]^(.6){\widetilde{\alpha_2}} & M_2 \ar[d]_{} \ar[r]^{} & \cdots \ar[d]_{} \ar[r]^{\alpha_{n-2}} &
  M_{n-2}
  \ar[d]_{} \ar[r]^{\alpha_{n-1}} & Y \ar[d]^{} \\
  \Sigma^{-1}F^iY \ar[r]^{} & F^iV \ar[r]^(.4){F^i\widetilde{\alpha_1}} & F^i(M_1\oplus M) \ar[r]^{} & F^iM_2
  \ar[r]^{} & \cdots
  \ar[r]^{} & F^iM_{n-2} \ar[r]^{F^i\alpha_{n-2}} & F^iY}$$

Let $p_X$ and $p_M$ be the projections of $V$ onto $X$ and $M$,
respectively. Since $\widetilde{\alpha_1}: V\ra M_1\oplus M$ is a
left $(\add(M), F, \Phi)$-approximation of $V$, $y_iF^ip_M$
factorizes through $\widetilde{\alpha_1}$. So there is a morphism
$s_i: M_1\oplus M\ra F^iM$ such that
$y_iF^ip_M=\widetilde{\alpha_1}s_i$. Hence
$x_0y_iF^ip_M=s(\Sigma^{-1}\widetilde{\alpha_n})\widetilde{\alpha_1}s_i=0$.
By assumption $X\in \Y(M)$, we have $\Hom_{\F}(M, F^iX)=0$. Then the
composition $vy_iF^ip_X: M'\stackrel{v}\ra V \stackrel{y_i}\ra F^iV
\stackrel{F^ip_X}\ra F^iX $ belongs to $\Hom_{\F}(M', F^iX)=0$, thus
$x_0y_iF^ip_X=uvy_iF^ip_X=0$. Altogether, we have shown that
$x_0y_i=0$ for $0\neq i\in \Phi$. Hence $xy\in I$, and $I$ is a
right ideal in $\E(V)$.

Obviously, $y_0x_0$ factorizes through an object in add$(M)$ and
through $\Sigma^{-1}\widetilde{\alpha_n}$. Set $0\neq i\in\Phi$.
Note that $\widetilde{\alpha_1}: V\ra M_1\oplus M$ is a left
$(\add(M), F, \Phi)$-approximation of $V$. Thus there is a morphism
$h_i: M_1\oplus M\ra F^iM'$ such that
$y_iF^iu=\widetilde{\alpha_1}h_i$. By assumption, we have
$\Hom_{\F}(M, F^iX)=0$ for $0\neq i\in \Phi$. This implies that
$h_iF^ivF^ip_X=0$, and therefore
$y_iF^ix_0F^ip_X=\widetilde{\alpha_1}h_iF^ivF^ip_X=0$. Since
$(\Sigma^{-1}\widetilde{\alpha_n})p_M=0$, we have shown that
$y_iF^ix_0F^ip_M=
y_iF^isF^i(\Sigma^{-1}\widetilde{\alpha_n})F^ip_M=y_iF^isF^i(\Sigma^{-1}\widetilde{\alpha_n}p_M)=0$.
Thus, $y_iF^ix_0=0$ for $0\neq i\in \Phi$. Hence $yx\in I$, and $I$
is a left ideal in $\E(V)$. Thus $I$ is an ideal in $\E(V)$.

In the same manner we can seen that $J$ is an ideal in $\E(W)$.
$\square$

\medskip
The following lemma is essentially taken from \cite{HKX}. The proof
remains valid for the present situation.

\begin{Lem}Then notations are the same as above. Then

$(1)$ $I\cdot\E(V,M)=0$.

$(2)$ $I\cdot\E(V,X)=\begin{array}{ll} \{(x_i)_{i\in \Phi}\in
\E(V,X)\mid x_i=0 \mbox{ for }0\neq i\in \Phi,& x_0 \mbox{
factorizes through }\\ &\add(M) \mbox{ and }\Sigma^{-1}\alpha_n\}
\end{array}$

$(3)$ For $x=(x_i)_{i\in \Phi}\in \E(V', X)$ with $V'\in$ $\add(V)$,
we have $\Imf(\mu(x))\subseteq I\cdot\E(V, X)$ if and only if
$x_i=0$ for all $0\neq i\in \Phi$ and $x_0$ factorizes through
$\add(M)$ and $\Sigma^{-1}\alpha_n$.

$(4)$ Let $f: M'\ra X$ with $M'\in \add(M)$. Then
$\Imf(\E(V,f))\subseteq I\cdot\E(V,X)$ if and only if $f$ factorizes
through $\Sigma^{-1}\alpha_n$.
  \label{lem1.2}
\end{Lem}

\medskip

Now, we turn to prove Theorem \ref{Theo}.

{\bf Proof of Theorem \ref{Theo}.} In order to prove Theorem
\ref{Theo}, Our strategy is trying to find out a tilting complex
over $\E(V)/I$ and compute its endomorphism ring. For convenience,
we define
$$\Lambda:=\E(V), \quad \Gamma:=\E(W), \quad \overline{\Lambda}:=\Lambda/I,\quad \overline{\Gamma}:=\Gamma/J.$$

Set
$$\begin{array}{ll}\widetilde{T^\bullet}: 0\ra \E(V,X)\stackrel{(V,\alpha_1)}\ra \E(V,M_1)\stackrel{(V,\alpha_2)}\ra
&\E(V,M_2)\stackrel{(V,\alpha_3)}\ra\cdots\\&\stackrel{(V,\overline{\alpha_{n-2}})}\ra\E(V,M_{n-2}\oplus
M)\ra 0.
\end{array}$$
Note that $\widetilde{T^\bullet}$ is a complex in
$K^b(\pmodcat{\Lambda})$. However, by easy computation,
$\widetilde{T^\bullet}$ is not a tilting complex over $\Lambda$.

Pick $x=(x_i)_{i\in \Phi}\in I\cdot \E(V,X)$. By the definition,
$\E(V,\alpha_1)(x)=(x_iF^i\alpha_1)_{i\in \Phi}$. Note that $x_i=0$
for $0\neq i\in \Phi$ and $x_0$ factorizes through
$\Sigma^{-1}\alpha_n$. So $\E(V, \alpha_1)(x)=0$. Hence the morphism
$\E(V,\alpha_1): \E(V, X)\ra \E(V, M_1)$ induces a morphism
$$q: \E(V,X)/I\cdot\E(V,X)\ra \E(V,M_1).$$

Let $P=\E(V,X)/I\cdot\E(V,X)$, and $p: \E(V,X)\ra
\E(V,X)/I\cdot\E(V,X)$ be the canonical surjective map. Then we can
write $\E(V,\alpha)=pq$.

Thus, we have a complex
$$T^\bullet: 0\ra P\stackrel{q}\ra \E(V,M_1)\stackrel{(V,\alpha_2)}\ra
\E(V,M_2)\stackrel{(V,\alpha_3)}\ra\cdots\stackrel{(V,\overline{\alpha_{n-2}})}\ra\E(V,M_{n-2}\oplus
M)\ra 0.$$ in $D^b(\overline{\Lambda})$. We will prove that
$T^\bullet$ is a tilting complex over $\overline{\Lambda}$.

Note that $\E(V,X)$ is a finitely generated projective left
$\Lambda$-module and $I\cdot\E(V,M)=0$. Then $P$ and $\E(V,M)$ are
finitely generated projective left $\overline{\Lambda}$-modules.
Hence $T^\bullet$ is a complex in
$K^b(\pmodcat{\overline{\Lambda}})$. Clearly, add$(T^\bullet)$
generates $K^b(\pmodcat{\Lambda})$ as a triangulated category. So it
suffices to prove that
$\Hom_{K^b(\overline{\Lambda})}(T^\bullet,T^\bullet[i])=0$ for
$i\neq 0$.

$(1)$ $\Hom_{K^b(\overline{\Lambda})}(T^\bullet, T^\bullet[i])=0$
for $i=1,2,\cdots,n-2$.

The first case: $i=1,2,\cdots,n-3$.

Let $f^\bullet$ be a morphism in
$\Hom_{K^b(\overline{\Lambda})}(T^\bullet, T^\bullet[i])$. For
simplicity, throughout the proof, we denote $\E(X,Y)$ by $(X,Y)$ in
commutative diagrams.
$$\xymatrix{
  0  \ar[r]^{} & P \ar[d]_{f^0} \ar[r]^{q} & (V,M_1) \ar@{-->}[ld]_{\mu(s^1)}\ar[d]_{f^1} \ar[r]^{(V,\alpha_2)}
  &(V,M_2)\ar@{-->}[ld]_{\mu(s^2)}\ar[d]_{f^2}\ar[r]^{(V,\alpha_3)}
  & \cdots  \ar[r]^{} &(V,M_{n-2-i})\ar@{-->}[ld]_{\mu(s^{n-2-i})}\ar[d]_{f^{n-2-i}}  \ar[r]^{} & \cdots  \\
  \cdots \ar[r]^{} &(V,M_i) \ar[r]^{(V,\alpha_{i+1})} & (V,M_{i+1}) \ar[r]^{} &(V,M_{i+2})\ar[r]
  &\cdots\ar[r]&(V,M_{n-2}\oplus M) \ar[r]^{} & 0   }$$

By Lemma $\ref{lem1.0}$$(1)$, we can assume that
$$\mu(x^0)=pf^0,\; f^{n-2-i}=\mu(x^{n-2-i}),\; f^j=\mu(x^j)$$
 with
$$x^0=(x^0_k)_{k\in \Phi}\in \E(X,M_i),$$
$$ \quad x^{n-2-i}=(x^{n-2-i}_k)_{k\in \Phi}\in
\E(M_{n-2-i},M_{n-2}\oplus M), $$
$$x^j=(x^j_k)_{k\in \Phi}\in \E(M_j,M_{i+j})$$
for $j=1,2,\cdots,n-3-i$.

Note that $\alpha_1: X\ra M_1$ is a left $(\add (M), F,
\Phi)$-approximation of $X$. Then there are morphisms $y^0_j: M_1\ra
F^jM_i$ such that $x^0_j=\alpha_1 y^0_j$ for $j\in \Phi$. We denote
$(y^0_j)_{j\in \Phi}$ by $y^0$.

Since
$$\E(V,\alpha_1)\mu(y^0)=\mu(\underline{\alpha_1})\mu(y^0)=\mu(\underline{\alpha_1}y^0)
=\mu((\alpha_1y^0_j)_{j\in \Phi})=\mu((x^0_j)_{j\in
\Phi})=\mu(x^0),$$ we can get $pq\mu(y^0)=\mu(x^0)=pf^0$. This
implies that $q\mu(y^0)=f^0$ since $p$ is surjective. We denote
$f^1-\mu(y^0)\E(V,\alpha_{i+1})$ by $s^1$.

Thus,
$$s^1=f^1-\mu(y^0)\E(V,\alpha_{i+1})=\mu(x^1)-\mu(y^0\underline{\alpha_{i+1}})=
\mu(x^1-y^0\underline{\alpha_{i+1}}).$$

We denote $x^1_j-y^0_jF^j\alpha_{i+1}$ by $s^1_j$. Note that
$$\E(V,\alpha_1)f^1=pqf^1=pf^0\E(V,\alpha_{i+1})=\mu(x^0)\E(V,\alpha_{i+1}),$$
i.e.,
$\mu(\underline{\alpha_1}x^1)=\mu(x^0\underline{\alpha_{i+1}})$.
This implies that $\alpha_1x^1_j=x^0_jF^j\alpha_{i+1}$ for $j\in
\Phi$.

It follows that
$$\alpha_1s^1_j=\alpha_1(x^1_j-y^0_jF^j\alpha_{i+1})=\alpha_1x^1_j-\alpha_1y^0_jF^j\alpha_{i+1}=
\alpha_1x^1_j-x^0_jF^j\alpha_{i+1}=0$$ for $j\in\Phi$. Then there
exists $y^1_j: M_2\ra M_{i+1}$ such that $s^1_j=\alpha_2y^1_j$ for
$j\in \Phi$. For convenience, we denote $(y^1_j)_{j\in \Phi}$ by
$y^1$. Now, we check that
$\E(V,\alpha_2)\mu(y^1)+\mu(y^0)\E(V,\alpha_{i+1})=f^1$.
$$\begin{array}{lll}
\E(V,\alpha_2)\mu(y^1)+\mu(y^0)\E(V,\alpha_{i+1})&=&\mu(\underline{\alpha_2})\mu(y^1)+\mu(y^0)\mu(\underline{\alpha_{i+1}})\\
&=&\mu(\underline{\alpha_2}y^1+y^0\underline{\alpha_{i+1}})\\
&=&\mu((\alpha_2y^1_j+y^0_jF^j\alpha_{i+1})_{j\in \Phi})\\
&=&\mu((x^1_j)_{j\in \Phi})\\
&=&f^1.
\end{array}$$

We denote $f^2-\mu(y^1)\E(V,\alpha_{i+2})$ by $s^2$. Thus,
$$s^2=f^2-\mu(y^1)\E(V,\alpha_{i+2})=\mu(x^2)-\mu(y^1)\mu(\underline{\alpha_{i+2}})=
\mu(x^2-y^1\underline{\alpha_{i+2}}).$$

We denote $x^2_j-y^1_jF^j\alpha_{i+2}$ by $s^2_j$ for $j\in \Phi$.
Note that $\E(V,\alpha_2)\mu(x^2)=\mu(x^1)\E(V,\alpha_{i+2})$, we
can get $\alpha_2x^2_j=x^1_jF^j(\alpha_{i+2})$ for $j\in \Phi$.

It follows that
$$\begin{array}{lll}
\alpha_2s^2_j&=&\alpha_2(x^2_j-y^1_jF^j\alpha_{i+2})\\
&=&\alpha_2x^2_j-s^1_jF^j\alpha_{i+2}\\
&=&\alpha_2x^2_j-(x^1_j-y^0_jF^j\alpha_{i+1})F^j\alpha_{i+2}\\
&=&\alpha_2x^2_j-x^1_jF^j\alpha_{i+2}\\
&=&0.
\end{array}$$

Hence there are $y^2_j: M_3\ra F^jM_{i+2}$ such that
$\alpha_3y^2_j=s^2_j$ for $j\in \Phi$. Similarly, we can check that
$f^2=\E(V,\alpha_3)\mu(y^2)+\mu(y^1)\E(V,\alpha_2)$. By induction,
we can prove that $f^\bullet$ is null-homotopic. Hence
$\Hom_{K^b(\overline{\Lambda})}(T^\bullet, T^\bullet[i])=0$ for
$i=1,2,\cdots,n-3$. The second case: $i=n-2$. It is easy to check
$\Hom_{K^b(\overline{\Lambda})}(T^\bullet, T^\bullet[n-2])=0$.

%

\medskip

$(2)$ $\Hom_{K^b(\overline{\Lambda})}(T^\bullet,T^\bullet[-i])=0$
for $i=1,\cdots,n-2$.

The first case: $i=1,\cdots,n-3.$ Let $f^\bullet$ be a morphism in
$\Hom_{K^b(\overline{\Lambda})}(T^\bullet, T^\bullet[i])$. We have
the following commutative diagram:
$$\xymatrix{
  \cdots  \ar[r]^{} & (V,M_i)\ar@{-->}[ldd]_(.3){\mu(g)} \ar[d]_{f^0} \ar[r]^{(V,\alpha_{i+1})} & (V,M_{i+1})\ar@{-->}[lldd]_(.3){\mu(s^1)} \ar[d]_{f^1} \ar[r]^{(V,\alpha_{i+2})}
  &(V,M_{i+2})\ar@{-->}[ld]_{\mu(s^2)}\ar[d]_{f^2}\ar[r]& \cdots  \ar[r]^{} & (V,M_{n-2}\oplus M)\ar@{-->}[ld]_{\mu(s^{n-i-2})}\ar[d]_{f^{n-i-2}}\ar[r]^{} & 0  \\
 0\ar[r]^{} & P\ar[r]^{q} & (V,M_1) \ar[r]^{(V,\alpha_2)} &(V,M_2) \ar[r]^{(V,\alpha_3)} & \cdots \ar[r]^{}
 & (V,M_{n-i-2})\ar[r]^{(V,\alpha_{n-i-1})}&\cdots\\
 (V, X)\ar[rru]_{(V, \alpha_1)}\ar[ur]^{p}&&&&&&}$$

By Lemma $\ref{lem1.0}$$(1)$, we assume $f^j=\mu(x^j),
f^{n-i-2}=\mu(x^{n-i-2})$ with $x^j=(x^j_k)_{k\in \Phi}\in
\E(M_{i+j},M_{j}),
x^{n-i-2}=(x^{n-i-2}_k)_{k\in\Phi}\in\E(M_{n-2}\oplus M, M_{n-i-2})$
for $j=1, 2, \cdots, n-i-3$. From the above commutative diagram, we
can get $f^{n-i-2}\E(V,\alpha_{n-i-1})=0$. This implies
$x^{n-i-2}_jF^j\alpha_{n-i-1}=0$ for $j\in \Phi$. So there are
morphisms
$$s^{n-i-2}_j: M_{n-2}\oplus M\ra F^j(M_{n-i-3})$$
such that
$$x^{n-i-2}_j=s^{n-i-2}_jF^j\alpha_{n-i-2}$$ for $j\in \Phi$.
We denote $(s^{n-i-2}_j)_{j\in \Phi}$ by $s^{n-i-2}$. So
$$\begin{array}{lll}
\mu(s^{n-i-2})\E(V,\alpha_{n-i-2})&=&\mu(s^{n-i-2}\underline{\alpha_{n-i-2}})\\
&=&\mu((s^{n-i-2}_jF^j\alpha_{n-i-2})_{j\in \Phi})\\
&=&\mu((x^{n-i-2}_j)_{j\in \Phi})\\
&=&f^{n-i-2}.
\end{array}$$

We denote $f^{n-i-3}-\E(V,\overline{\alpha_{n-2}})\mu(s^{n-i-2})$ by
$t^{n-i-3}$, and we write $t^{n-i-3}_j$ instead of
$x^{n-i-3}_j-\overline{\alpha_{n-2}}s^{n-i-2}_j$ for $j\in \Phi$.
Note that
$$\E(V,\overline{\alpha_{n-2}})f^{n-i-2}=f^{n-i-3}\E(V,\alpha_{n-i-2}).$$
Then
$$(\overline{\alpha_{n-2}}x^{n-i-2}_j)_{j\in
\Phi}=(x^{n-i-3}_jF^j\alpha_{n-i-2})_{j\in \Phi}.$$

We can deduce
$$\begin{array}{lll}
t^{n-i-3}_jF^j\alpha_{n-i-2}&=&
(x^{n-i-3}_j-\overline{\alpha_{n-2}}s^{n-i-2}_j)F^j\alpha_{n-i-2}\\
&=&x^{n-i-3}_jF^j\alpha_{n-i-2}-\overline{\alpha_{n-2}}x^{n-i-2}_j\\
&=&0.
\end{array}$$

So there exist morphisms $s^{n-i-3}_j: M_{n-3}\ra F^jM_{n-i-4}$ such
that $$s^{n-i-3}_jF^j\alpha_{n-i-3}=t^{n-i-3}_j$$ for $j\in\Phi$.

We denote $(s^{n-i-3}_j)_{j\in \Phi}$ by $s^{n-i-3}$. We can deduce
$$\begin{array}{ll}
&\mu(s^{n-i-3})\E(V,\alpha_{n-i-3})+\E(V,\overline{\alpha_{n-2}})\mu(s^{n-i-2})\\=&
\mu(s^{n-i-3})\mu(\underline{\alpha_{n-i-3}})+\mu(\underline{\overline{\alpha_{n-2}}})
\mu(s^{n-i-2})\\
=&\mu((s^{n-i-3}_jF^j\alpha_{n-i-3}+\overline{\alpha_{n-2}}s^{n-i-2}_j)_{j\in
\Phi})\\
=&\mu((x^{n-i-3}_j)_{j\in \Phi})\\
=&f^{n-i-3}.
\end{array}$$

By induction, there are morphisms
$$\mu(s^1): \E(V, M_{i+1})\ra \E(V, X),\; \mu(s^k): \E(V,M_{i+k})\ra \E(V,M_k-1)$$ such that
$$f^1=\mu(s^1)\E(V, \alpha_1)+\E(V, \alpha_{i+2})\mu(s^2), \;
f^k=\mu(s^k)\E(V,\alpha_k)+\E(V,\alpha_{i+k+1})\mu(s^{k+1})$$ for
$k=2,\cdots,n-i-2$. Here we define $\mu(s^{n-i-1})=0$.

Hence it suffices to prove that $f^0=\E(V,\alpha_{i+1})\mu(s^1)p$.
Note that $p: \E(V,X)\ra P$ is surjective and $\E(V,M_i)$ is a
projective $\E(V)$-module. Then there exists a morphism $\mu(g):
\E(V,M_i)\ra \E(V,X)$ such that $f^0=\mu(g)p$. Since $X\in \Y(M)$,
we can get $g_j=0, s^1_j=0$ for $0\neq j\in \Phi$. Note that
$f^0q=\E(V,\alpha_{i+1})f^1$, this implies
$g_0\alpha_1=\alpha_{i+1}x^1_0$. So
$$(\alpha_{i+1}s^1_0-g_0)\alpha_1=\alpha_{i+1}(x^1_0-\alpha_{i+2}s^2_0)-g_0\alpha_1
=\alpha_{i+1}x^1_0-g_0\alpha_1=0.$$ This implies that
$\alpha_{i+1}s^1_0-g_0$ can factorizes through
$\Sigma^{-1}\alpha_n$. By Lemma $\ref{lem1.2}$$(4)$, we can get
$\Imf \mu(\underline{\alpha_{i+1}s^1_0-g_0})\subseteq
I\cdot\E(V,X)$. So $(\mu(g)-\E(V,\alpha_{i+1})\mu(s^1))p=0$. This
implies $f^0=\E(V, \alpha_{i+1})\mu(s^1)p$. Hence $f^\bullet$ is
null-homotopic. The second case: $i=n-2$. We can verify similarly.

Hence $T^\bullet$ is a tilting complex over $\overline{\Lambda}$.


\medskip

%


Clearly, the homotopy category $K^b(\overline{\Lambda})$ can be
viewed as a full subcategory of $K^b(\Lambda)$. Thus, we have a ring
isomorphism
$\End_{K^b(\Pmodcat{\overline{\Lambda}})}(T^\bullet)\simeq
\End_{K^b(\Lambda)}(T^\bullet)$. Now, we will determine the
endomorphism ring $\End_{K^b(\Lambda)}(T^\bullet)$.

Let $f^\bullet\in \End_{K^b(\Lambda)}(T^\bullet)$. There is an
$\Lambda$-homomorphism $u^0: \E(V,X)\ra \E(V,X)$ such that
$u^0p=pf^0$, because $p: \E(V,X)\ra P$ is an epimorphism and
$\E(V,X)$ is a projective $\Lambda$-module. By Lemma
$\ref{lem1.0}$$(1)$, we can assume
$$u^0=\mu(x^0), f^{n-2}=\mu(x^{n-2}), f^i=\mu(x^i)$$
with
$$x^0=(x^0_i)_{i\in \Phi}\in \E(X), x^i=(x^i_j)_{j\in
\Phi}\in \E(M_i,M_i),$$
$$x^{n-2}=(x^{n-2}_j)_{j\in \Phi}\in
\E(M_{n-2}\oplus M)$$ for $i=1,\cdots,n-3$.

$\xymatrix@R=2mm{(V, X)\ar[dr]_(.6){p}\ar[drr]^{(V,
\alpha_1)}\ar@{-->}[dd]_(.3){u^0=\mu(x^0)}&&&&&&\\
\ar[r]&P\ar[r]_(.3){q}\ar[dd]_(.3){f^0}&(V, M_1)\ar[r]^{(V,
\alpha_2)}\ar[dd]^{f^2=\mu(x^1)}&(V,
M_2)\ar[r]\ar[dd]^{f^2=\mu(x^2)}&\cdots\ar[r]^(.2){(V,
\overline{\alpha_{n-2}})}&(V, M_{n-2}\oplus
M)\ar[dd]^{f^{n-2}=\mu(x^{n-2})}\ar[r]&0\\
(V, X)\ar[dr]_(.6){p}\ar[drr]^{(V, \alpha_1)}&&&&&&\\
\ar[r]&P\ar[r]_(.3){q}&(V, M_1)\ar[r]^{(V, \alpha_2)}&(V,
M_2)\ar[r]&\cdots\ar[r]^(.2){(V, \overline{\alpha_{n-2}})}&(V,
M_{n-2}\oplus M)\ar[r]&0 }$

By the commutativity of the above diagram, we have
$$\E(V,\alpha_1)f^1=\mu(x^0)\E(V,\alpha_1),$$
$$\E(V,\alpha_i)f^i=f^{i-1}\E(V,\alpha_i) \mbox{ for }i=2,\cdots,n-3,$$
$$\E(V,\overline{\alpha_{n-2}})f^{n-2}=f^{n-3}\E(V,\overline{\alpha_{n-2}}).$$
It follows that
$$\alpha_1x^1_j=x^0_jF^j\alpha_1,$$
$$ \alpha_ix^i_j=x^{i-1}_jF^j\alpha_i \mbox{ for }i=2,\cdots,n-3,$$
$$\overline{\alpha_{n-2}}x^{n-2}_j=x^{n-3}_jF^j\overline{\alpha_{n-2}} $$
for $j\in \Phi$ from Lemma $\ref{lem1.0}$$(1)$.

By Lemma $\ref{lem1.7}$, we can form the following commutative
diagram in $\F$:
$$\xymatrix{
  X \ar[d]_{x^0_i} \ar[r]^{\alpha_1} & M_1 \ar[d]_{x^1_i} \ar[r]^{\alpha_2} &
  M_2 \ar[d]_{x^2_i} \ar[r]^{\alpha_3} &
  \cdots  \ar[r]^(.3){\overline{\alpha_{n-2}}} &
  M_{n-2}\oplus M  \ar[d]_{x^{n-2}_i} \ar[r]^(.6){\overline{\alpha_{n-1}}} &
  W \ar@{-->}[d]_{h_i} \ar[r]^{\overline{\alpha_n}} & \Sigma X \ar[d]^{\Sigma x^0_i} \\
  F^iX \ar[r]^{F^i\alpha_1} & F^iM_1 \ar[r]^{F^i\alpha_2} &
  F^iM_2 \ar[r]^{F^i\alpha_3} &
  \cdots \ar[r]^(.2){F^i(\overline{\alpha_{n-2}})} &
  F^i(M_{n-2}\oplus M) \ar[r]^(.7){F^i\overline{\alpha_{n-1}}} & F^iW \ar[r]^{} & \Sigma(F^iX)
  }\quad\quad(\star)$$
where $h_i\in \Hom_{\F}(W,F^iW)$. Thus, for each $f^\bullet\in
\End_{K^b(\Lambda)}(T^\bullet)$, we can get an element
$h:=(h_i)_{i\in \Phi}\in \Gamma$. Define the following
correspondence:
$$\Theta: \End_{K^b(\Lambda)}(T^\bullet)\ra \overline{\Gamma}=\Gamma/J,$$
$$f^\bullet\mapsto h+J.$$

Now, we will prove that the correspondence $\Theta$ is a ring
homomorphism. The proof is divided into four steps.

Step 1. we will prove that $\Theta$ is well-defined. Suppose that
$f^\bullet\in \End_{K^b(\Lambda)}(T^\bullet)$ is null-homotopic,
that is, there are
$$ r_1: \E(V,M_1)\ra P, \; r_i: \E(V,M_{i})\ra \E(V,M_{i-1}), i=2,\cdots,n-3,$$
$$r_{n-2}: \E(V,M_{n-2}\oplus M)\ra\E(V,M_{n-3}),$$\\
such that
$$f^0=qr_1, \; f^1=r_1q+\E(V,\alpha_2)r_2,$$
$$f^i=r_i\E(V,\alpha_i)+\E(V,\alpha_{i+1}r_{i+1})r_{i+1}\mbox{   for } i=2,\cdots,n-3,$$
$$f^{n-3}=r_{n-3}\E(V,\alpha_{n-3})+\E(V,\overline{\alpha_{n-2}})r_{n-2}, \; f^{n-2}=r_{n-2}\E(V,\overline{\alpha_{n-2}}).$$

Since $p$ is surjective and $\E(V,M_1)$ is projective, there is a
morphism $s: \E(V,M_1)\ra \E(V,X)$ such that $r_1=sp$. By Lemma
$\ref{lem1.0}$$(1)$, we can assume
$$s=\mu(t), r_{n-2}=\mu(l)$$
with
$$t=(t_i)_{i\in \Phi}\in \E(M_1,X), l=(l_i)_{i\in \Phi}\in \E(M_{n-2}\oplus M, M_{n-3}).$$

$\xymatrix@R=2mm{(V, X)\ar[dr]_(.6){p}\ar[dd]_(.3){u^0=\mu(x^0)}&&&&&&\\
\ar[r]&P\ar[r]^(.3){q}\ar[dd]_{f^0}&(V,
M_1)\ar@{-->}[dll]_(.6){s=\mu(t)}\ar[ddl]^{r_1}\ar[r]^{(V,
\alpha_2)}\ar[dd]^{f^2=\mu(x^1)}&(V,
M_2)\ar[ddl]^{r_2}\ar[r]\ar[dd]^{f^2=\mu(x^2)}&\cdots\ar[r]^(.2){(V,
\overline{\alpha_{n-2}})}&(V, M_{n-2}\oplus
M)\ar[ddl]^{r_{n-2}=\mu(l)}\ar[dd]^{f^{n-2}=\mu(x^{n-2})}\ar[r]&0\\
(V, X)\ar[dr]_(.6){p}&&&&&&\\
\ar[r]&P\ar[r]_(.3){q}&(V, M_1)\ar[r]_{(V, \alpha_2)}&(V,
M_2)\ar[r]&\cdots\ar[r]_(.2){(V, \overline{\alpha_{n-2}})}&(V,
M_{n-2}\oplus M)\ar[r]&0 }$

By the definition of $\Y(M)$, we have $t_i=0$ for $0\neq i\in \Phi$.
It follows that
$$\mu(x^0-\underline{\alpha_1t_0})p=(u^0-pqs)p=0, \; \mu(x^{n-2})=\mu(l)\E(V,\overline{\alpha_{n-2}}).$$

It follows immediately that
$$\Imf \mu(x^0-\underline{\alpha_1t_0})\subseteq I\cdot\E(V,X), \; (x^{n-2}_i)_{i\in \Phi}=(l_iF^i\overline{\alpha_{n-2}})_{i\in
\Phi}.$$

By Lemma $\ref{lem1.2}$$(2)$, we can get that $x^0_i=0$ for $0\neq
i\in \Phi$ and $x^0_0-\alpha_1t_0$ factorizes through add$(M)$ and
$\Sigma^{-1}\alpha_n$. So $x_0^0-\alpha_1t_0=ab$ for some morphisms
$a: X\ra M'$ and $b: M'\ra X$ with $M'\in \add(M)$. Since $\alpha_1:
X\ra M_1$ is a left $(\add(M),F,\Phi)$-approximation of $X$, there
is a morphism $c: M_1\ra M'$ such that $a=\alpha_1c$. It follows
that
$$x^0_0=ab+\alpha_1t_0=\alpha_1cb+\alpha_1t_0=\alpha_1(cb+t_0).$$

Since
$\overline{\alpha_{n-1}}h_i=x^{n-2}_iF^i\overline{\alpha_{n-1}}=
l_iF^i\overline{\alpha_{n-2}}F^i\overline{\alpha_{n-1}}=0$, $h_i$
factorizes through $\overline{\alpha_n}$. So $h_i|_M=0$ since
$\overline{\alpha_n}|_M=0$. Since $x^0_i=0$ for $0\neq i\in \Phi$
and $Y\in \X(M)$, we deduce $h_i|_Y=0$. It follows that $h_i=0$ for
$0\neq i\in \Phi$.

We have
$\overline{\alpha_{n-1}}h_0=x^{n-2}_0\overline{\alpha_{n-1}}=
l_0\overline{\alpha_{n-2}}\overline{\alpha_{n-1}}=0$ which implies
that $h_0$ factorizes through $\overline{\alpha_n}$. Since
$h_0\overline{\alpha_n}=\overline{\alpha_n}\Sigma
x^0_0=\overline{\alpha_n}(\Sigma\alpha_1)\Sigma(cb+t)=0$, the
morphism $h_0$ factorizes through $M_{n-2}\oplus M$ which is in
add$(M)$. Thus, $h$ is an element in $J$. So $\Theta$ is
well-defined.


Step 2. we will prove that the map $\Theta$ is injective. Suppose
that $\Theta(f^\bullet)=h+J=J$. It suffices to prove that
$f^\bullet$ is null-homotopic. By the definition of $J$, we have
that $h_i=0$ for $0\neq i\in \Phi$, and $h_0$ factorizes through
add$(M)$ and $\overline{\alpha_n}$. Since $h_i=0$ for $0\neq i\in
\Phi$ and $h_0$ factorizes through $\overline{\alpha_n}$, we have
$x^{n-2}_iF^i\overline{\alpha_{n-1}}=0$, by the commutativity of
$(\star)$. Thus, there is a morphism $r^{n-2}_i: M_{n-2}\oplus M\ra
F^iM_{n-3}$ such that
$x^{n-2}_i=r^{n-2}_iF^i\overline{\alpha_{n-2}}$ for $i\in\Phi$. Let
us denote $r^{n-2}$ the morphism $(r^{n-2}_i)_{i\in \Phi}$. Then
$\mu(r^{n-2})\E(V,\overline{\alpha_{n-2}})=\mu(x^{n-2})$. And we
will denote $s^{n-3}$ the morphism
$f^{n-3}-\E(V,\overline{\alpha_{n-2}})\mu(r^{n-2})$. Thus
$s^{n-3}_i=x^{n-3}_i-\overline{\alpha_{n-2}}r^{n-2}_i$ for
$i\in\Phi$. Since
$f^{n-3}\E(V,\overline{\alpha_{n-2}})=\E(V,\overline{\alpha_{n-2}})f^{n-2}$,
that is, $(x^{n-3}_iF^i\overline{\alpha_{n-2}})_{i\in \Phi}=
(\overline{\alpha_{n-2}}x^{n-2}_i)_{i\in \Phi}$, we can deduce
$$\begin{array}{lll}
(x^{n-3}_i-\overline{\alpha_{n-2}}r^{n-2}_i)F^i(\overline{\alpha_{n-2}})&
=&x^{n-3}_iF^i\overline{\alpha_{n-2}}-\overline{\alpha_{n-2}}r^{n-2}_iF^i\overline{\alpha_{n-2}}\\
&=&x^{n-3}_iF^i\overline{\alpha_{n-2}}-\overline{\alpha_{n-2}}x^{n-2}_i\\
&=&0.
\end{array}$$

Thus, there are morphisms $r^{n-3}_i: M_{n-3}\ra F^iM_{n-4}$ such
that
$x^{n-3}_i-\overline{\alpha_{n-2}}r^{n-2}_i=r^{n-3}_iF^i\alpha_{n-3}$
for $i\in \Phi$. We denote $(r^{n-3}_i)_{i\in \Phi}$ by $r^{n-3}$.

Note that
$x^{n-3}_i-\overline{\alpha_{n-2}}r^{n-2}_i=r^{n-3}_iF^i\alpha_{n-3}$
for $i\in \Phi$. Then
$$\begin{array}{lll}
\mu(r^{n-3})\E(V,\alpha_{n-3})+\E(V,\overline{\alpha_{n-2}})\mu(r^{n-2})&=
&\mu((r^{n-3}_iF^i\alpha_{n-3}+\overline{\alpha_{n-2}}r^{n-2}_i)_{i\in
\Phi})\\
&=&\mu((x^{n-3}_i)_{i\in \Phi})\\
&=&\mu(x^{n-3})\\
&=&f^{n-3}.
\end{array}$$

By induction, we can construct
$$r^i:=(r^i_j)_{j\in \Phi}\in \E(M_i,M_{i-1})$$
and
$$s^i:=f^i-\E(V,\alpha_{i+1})\mu(s^{i+1})=(f^i_j-\alpha_{i+1}s^{i+1}_j)_{j\in
\Phi}\in \E(M_i,M_i)$$ satisfying that
$$f^i=\mu(s^i)\E(V,\alpha_{i})+\E(V,\alpha_{i+1})\mu(s^{i+1})$$
for $i=2,\cdots,n-4$. Let us denote $s^1$ the morphism
$$f^1-\E(V,\alpha_2)\mu(r^2)=(f^1_i-\alpha_2r^2_i)_{i\in \Phi}\in
\E(M_1,M_1).$$

Note that $\E(V,\alpha_2)f^2=f^1\E(V,\alpha_2)$, that is
$(\alpha_2x^2_i)_{i\in \Phi}=(x^1_iF^i\alpha_2)_{i\in \Phi}$. Then
$$\begin{array}{lll}
s^1_iF^i\alpha_2&=&(x^1_i-\alpha_2r^2_i)F^i\alpha_2\\
&=&x^1_iF^i\alpha_2-\alpha_2r^2_iF^i\alpha_2\\
&=&x^1_iF^i\alpha_2-\alpha_2(x^2_i-\alpha_3r^3_i) \\
&=&x^1_iF^i\alpha_2-\alpha_2x^2_i\\
&=&0.
\end{array}$$

Thus, there are morphisms $r^1_i: M_1\ra F^iX$ such that
$r^1_iF^i\alpha_1=s^1_i=x^1_i-\alpha_2r^2_i$ for $i\in \Phi$. We
define $r^1:=(r^1_i)_{i\in \Phi}$. Since $X\in \Y(M)$, we have
$r^1_i=0$ for $0\neq i\in \Phi$. Consequently,
$$f^1=\E(V,\alpha_2)\mu(r^2)+\mu(r^1)\E(V,\alpha_1).$$

We can get $\overline{\alpha_n}\Sigma x^0_i=0$ by the assumption
that $h_i=0$ for $0\neq i\in \Phi$. Thus, $x^0_i$ factorizes through
$\alpha_1$. Since $X\in \Y(M)$, we can obtain $x^0_i=0$ for $0\neq
i\in \Phi$.

Note that $u^0\E(V,\alpha_1)=\E(V,\alpha_1)f^1$. Then
$$\begin{array}{lll}
(x^0_0-\alpha_1r^1_0)\alpha_1&=&x^0_0\alpha_1-\alpha_1r^1_0\alpha_1\\
&=&x^0_0\alpha_1-\alpha_1(x^1_0-\alpha_2r^2_0)\\
&=&x^0_0\alpha_1-\alpha_1x^1_0\\
&=&0.
\end{array}$$
This implies that $x^0_0-\alpha_1r^1_0$ factorizes through
$\Sigma^{-1}\overline{\alpha_n}$.

Now, we prove that $x^0_0-\alpha_1r^1_0$ factorizes through
add$(M)$. Since $\alpha_1r^1_0$ factorizes through add$(M)$, it
suffices to prove that $x^0_0$ can factorize through add$(M)$. By
assumption, $h_0$ can factorize through add$(M)$. So there are
morphisms $a: W\ra M'$ and $b: M'\ra W$ such that $h_0=ab$ for
$M'\in$ add$(M)$. Since $\overline{\alpha_{n-1}}$ is right
$(\add(M), F,-\Phi)$-approximation of $W$, there is a morphism $c:
M'\ra M_{n-2}\oplus M$ such that $b=c\overline{\alpha_{n-1}}$.
Consequently,
$$\overline{\alpha_n} \Sigma
x^0_0=h_0\overline{\alpha_n}=ac\overline{\alpha_{n-1}}\overline{\alpha_n}=0.$$
This implies that $x^0_0$ can factorize through $M_1$ which belongs
to add$(M)$. By Lemma $\ref{lem1.2}$$(3)$, we deduce
$\Imf\mu(x^0-\underline{\alpha_1r^1_0})\subseteq I\cdot \E(V,X)$.
Thus,
$$p(f^0-q\mu(\underline{r^1_0})p)=pf^0-pq\mu(\underline{r^1_0})p=
(u^0-pq\mu(\underline{r^1_0}))p=0.$$

Hence $f^0=q\mu(\underline{r^1_0})p$. Altogether, we have proven
that $f^\bullet$ is null-homotopic.

Step 3. we will prove that the map $\Theta$ is surjective. Let
$h=(h_i)_{i\in \Phi}$ with $h_i: W\ra F^iW$ for $i\in \Phi$. Since
$\overline{\alpha_{n-1}}$ is a right
$(\add(M),F,-\Phi)$-approximation of $W$, there is a commutative
diagram:
$$\xymatrix{
  X \ar[d]_{x^0_i} \ar[r]^{\alpha_1} & M_1 \ar[d]_{x^1_i} \ar[r]^{\alpha_2} &
  \cdots\ar[r]^{\alpha_{n-3}} &
  M_{n-3} \ar[d]_{x^{n-3}_i} \ar[r]^(.4){\overline{\alpha_{n-2}}}
  &M_{n-2}\oplus M \ar[d]_{x^{n-2}_i} \ar[r]^(.6){\overline{\alpha_{n-1}}}
  & W \ar[d]_{h_i} \ar[r]^{\overline{\alpha_n}} & \Sigma X \ar[d]^{\Sigma x^0_i} \\
  F^iX \ar[r]^{F^i\alpha_1} & F^iM_1 \ar[r]^{F^i\alpha_2} & \cdots \ar[r]^{} & F^iM_{n-3} \ar[r]^(.4){F^i\overline{\alpha_{n-2}}}
   & F^i(M_{n-2}\oplus M) \ar[r]^(.6){F^i\overline{\alpha_{n-1}}}
   & F^iW \ar[r]^{} & \Sigma F^iX   }$$

We denote $(x^j_i)_{i\in \Phi}$ by $x^j$ for $j=0,1,\cdots,n-2$.
From the commutative diagram, we have
$\overline{\alpha_{n-2}}x^{n-2}_i=x^{n-3}_iF^i\overline{\alpha_{n-2}}$
and $\alpha_jx^j_i=x^{j-1}_iF^i\alpha_j$ for $j=1,2,\cdots,n-2$.
This implies
$$\E(V,\overline{\alpha_{n-2}})\mu(x^{n-2})=\mu(x^{n-3})\E(V,\overline{\alpha_{n-2}}),$$
$$\mu(\underline{\alpha_j})\mu(x^j)=\mu(x^{j-1})\mu(\underline{\alpha_j}) \mbox{ for }j=1,\cdots,n-2.$$

So we have the following commutative diagram
$$\xymatrix@R=2mm{(V, X)\ar[dr]_(.6){p}\ar[drr]^{(V,
\alpha_1)}\ar[dd]_(.3){u^0=\mu(x^0)}&&&&&&\\
\ar[r]&P\ar[r]_(.3){q}\ar@{-->}[dd]_(.3){f^0}&(V, M_1)\ar[r]^{(V,
\alpha_2)}\ar[dd]^{f^2=\mu(x^1)}&(V,
M_2)\ar[r]\ar[dd]^{f^2=\mu(x^2)}&\cdots\ar[r]^(.2){(V,
\overline{\alpha_{n-2}})}&(V, M_{n-2}\oplus
M)\ar[dd]^{f^{n-2}=\mu(x^{n-2})}\ar[r]&0\\
(V, X)\ar[dr]_(.6){p}\ar[drr]^{(V, \alpha_1)}&&&&&&\\
\ar[r]&P\ar[r]_(.3){q}&(V, M_1)\ar[r]^{(V, \alpha_2)}&(V,
M_2)\ar[r]&\cdots\ar[r]^(.2){(V, \overline{\alpha_{n-2}})}&(V,
M_{n-2}\oplus M)\ar[r]&0 }$$

We conclude from $\mu(x^0)(I\cdot\E(V,X))\subseteq I\cdot\E(V,X)$,
that $\mu(x^0)$ induces a morphism $f^0: P\ra P$ satisfying that
$pf^0=\mu(x^0)p$, and finally that
$$p(f^0q-q\mu(x^1))=\mu(x^0)\E(V,\alpha_1)-\E(V,\alpha_1)\mu(x^1)=0.$$

Note that $p$ is surjective. Then $f^0q=q\mu(x^1)$. Define
$f^i=\mu(x^i)$ for $i=1,\cdots,n-2$. Hence $\Theta$ is surjective.

Step 4. we will prove that the map $\Theta$ is a ring homomorphism.
Take $f^\bullet$ and $g^\bullet$ in
$\End_{K^b(\Lambda)}(T^\bullet)$. Since $p$ is surjective and
$\E(V,X)$ is projective as left $\E(V)$-module, there is a map
$\mu(x^0): \E(V,X)\ra \E(V,X)$ such that $\mu(x^0)p=pf^0$.
Similarly, there is a map $\mu(y^0): \E(V,X)\ra \E(V,X)$ such that
$\mu(y^0)p=pg^0$. Suppose that $f^i=\mu(x^i), g^i=\mu(y^i)$ for
$i=1,\cdots,n-2$. Define $h:=(h_i)_{i\in \Phi}$ and
$h':=(h'_i)_{i\in \Phi}$ be in $\Gamma$ such that
$$\begin{array}{ll}
\overline{\alpha_{n-1}}h_i=x^{n-2}_iF^i\overline{\alpha_{n-1}},
  & \overline{\alpha_n}\Sigma x^0_i=h_i(F^i\overline{\alpha_n})\delta(F,i,X,1)  \\
  \overline{\alpha_{n-1}}{h'}_i=y^{n-2}_iF^i\overline{\alpha_{n-1}},
    & \overline{\alpha_n}\Sigma y^0_i={h'}_i(F^i\overline{\alpha_n})\delta(F,i,X,1)
\end{array}$$
for $i\in \Phi$. By definition, we have $\Theta(f^\bullet)=h+J,
\Theta(g^\bullet)=h'+J$ and
$$\Theta(f^\bullet)\Theta(g^\bullet)=(\sum_{i,j\in \Phi \atop i+j=k} h_iF^ih'_j)_{k\in \Phi}+J.$$

Now, we calculate $\Theta(f^\bullet g^\bullet)$.
$$x^{n-2}y^{n-2}=(\sum_{i,j\in \Phi \atop
i+j=k}x^{n-2}_iF^iy^{n-2}_j)_{k\in \Phi}, \; x^0y^0=(\sum_{i,j\in
\Phi \atop i+j=k}x^0_iF^iy^0_j)_{k\in \Phi }.$$

For each $k\in \Phi$
$$\begin{array}{lll}
 \overline{\alpha_{n-1}}(\sum_{i,j\in \Phi \atop i+j=k}
 h_iF^ih'_j)&=&\sum_{i,j\in \Phi \atop
 i+j=k}\overline{\alpha_{n-1}}h_iF^ih'_j\\
 &=&\sum_{i,j\in \Phi\atop
 i+j=k}x^{n-2}_iF^i(y^{n-2}_jF^j\overline{\alpha_{n-1}})\\
 &=&\sum_{i,j\in \Phi\atop
 i+j=k}x^{n-2}_iF^iy^{n-2}_jF^k\overline{\alpha_{n-1}.}
\end{array}$$
$$\begin{array}{lll}
  \overline{\alpha_n}\Sigma (x^0y^0)_k&=&
  \overline{\alpha_n}\Sigma(\sum_{i,j\in \Phi \atop
  i+j=k}x^0_iF^iy^0_j)\\
&=&\sum_{i,j\in \Phi \atop i+j=k}\overline{\alpha_n}(\Sigma
x^0_i)(\Sigma F^iy^0_j)\\
&=&\sum_{i,j\in \Phi \atop
i+j=k}h_iF^i\overline{\alpha_n}\delta(F,i,X,1)(\Sigma F^iy^0_j)\\
&=&\sum_{i,j\in \Phi \atop
i+j=k}h_iF^j{h'}F^k\overline{\alpha_n}\delta(F,k,X,1).
\end{array}$$
So $\Theta(f^\bullet g^\bullet)=\Theta(f^\bullet)\Theta(g^\bullet)$.
Thus, $\Theta$ is a ring homomorphism. $\square$

\medskip

\medskip

If $\F$ is a triangulated $R$-category, we can get the main result
in \cite{HKX}. Combined with \cite[Theorem 3.1]{GKO}, we can get the
following corollary.




\smallskip


\begin{Koro} \label{coro2}Let $\Phi$ be an admissible subset of $\mathbb{N}$. Let $\mathcal {F}_3$ be a triangulated $k$-category with an
$(n-2)$-cluster tilting subcategory $\F$, which is closed under
$\Sigma_3^{n-2}$, where $\Sigma_3$ denotes the suspension functor in
$\mathcal {F}_3$. Suppose that there exists a diagram
$$\xymatrix{
   & X_2 \ar[dr]\ar[rr]^{\alpha_2} && X_3 \ar[dr] \ar[rr] && X_4 \ar[dr]\ar[r]& \cdots&\cdots\ar[r]&X_{n-1} \ar[dr]^{\alpha_{n-1}}& \\
 X_1\ar[ur]^{\alpha_1}&&X_{2.5}\ar[ur] \ar@{-->}[ll]&&X_{3.5}\ar[ur] \ar@{-->}[ll]&&\cdots\ar@{-->}[ll]
 &X_{n-1.5}\ar@{-->}[l] \ar[ur]&&X_n\ar@{-->}[ll]}$$
in $\mathcal {F}_3$, satisfying that

$(1)$ $\alpha_1: X_1\ra X_2$ is a left $(\add (X), F,
\Phi)$-approximation of $X_1$

$(2)$ $\alpha_{n-1}: X_{n-1}\ra X_n$ is a right $(\add (X), F,
-\Phi)$-approximation of $X_n$,

$(3)$ $X_1\in \mathscr{Y}^{\Sigma_{3}^{n-2},\Phi}_{\F}(X)$,
$X_{n-1}\in \mathscr{X}^{\Sigma_{3}^{n-2},\Phi}_{\F}(X)$,

\noindent where $X$ is the direct sum of $X_i$ for
$i=2,3,\cdots,n-1$.

Then we can get that the two algebras ${\rm
E}^{\Sigma_{3}^{n-2},\Phi}_{\F}(X_1\oplus X)/I$ and ${\rm
E}^{\Sigma_{3}^{n-2},\Phi}_{\F}(X_{n-1}\oplus X)/J$ are derived
equivalent, where $\mathscr{X}^{\Sigma_{3}^{n-2},\Phi}_{\F}(X),
\mathscr{Y}^{\Sigma_{3}^{n-2},\Phi}_{\F}(X), I$ and $J$ are defined
as in Theorem $\ref{Theo}$.
\end{Koro}

{\bf Proof.} This follows from \cite[Theorem 3.1]{GKO} and Theorem
\ref{Theo}. $\square$


\medskip
In \cite{IY}, Iyama and Yoshino introduced Auslander-Reiten
$n$-angles in $(n-2)$-cluster tilting subcategories of triangulated
$k$-categories and proved that they always exist. Let $\mathcal {T}$
be a Krull-Schmidt triangulated category with shift functor
$\Sigma_3$, and let $\mathcal {S}$ be an $n$-cluster tilting
subcategory of $\mathcal {T}$.
$$X_{i+1}\stackrel{b_{i+1}}\ra C_i\stackrel{a_i}\ra X_i\ra \Sigma_3X_{i+1}\quad (0\leq i<n).$$
are triangles in $\mathcal {T}$. A complex
$$X_n\stackrel{b_n}\ra C_{n-1}\stackrel{a_{n-1}b_{n-1}}\ra C_{n-2}\stackrel{a_{n-2}b_{n-2}}\ra\cdots\stackrel{a_2b_2}\ra C_1\stackrel{a_1b_1}
\ra C_0\stackrel{a_0}\ra X_0$$ is called an {\em Auslander-Reiten
$(n+2)$-angle} if the following conditions are satisfied.

\quad$(1)$ $X_n, X_0$ and $C_i (0\leq i< n)$ belong to $\mathcal
{S}$.

\quad$(2)$ $a_0$ is a sink map of $X_0$ in $\mathcal{S}$ and $b_n$
is a source map of $X_n$ in $\mathcal{S}$.

\quad$(3)$ $a_i$ is a minimal right $\mathcal {S}$-approximation of
$X_i$ for $0<i<n$.

\quad$(4)$ $b_i$ is a minimal left $\mathcal{S}$-approximation of
$X_i$ for $0<i<n$.

%

As a corollary of Corollary \ref{coro2}, we can establish a
relationship between Auslander-Reiten $n$-angle and derived
equivalences.

\begin{Koro}\label{nAa}
Let $\mathcal{T}$ be a Krull- Schmidt triangulated $k$-category with
shift functor $\Sigma_3$, and let $\mathcal {S}$ be an
$(n-2)$-cluster tilting subcategory of $\mathcal {T}$, which is
closed under $\Sigma^{n-2}_3$. Suppose that
$$X_1\stackrel{\alpha_1}\ra X_2\stackrel{\alpha_2}\ra X_3\ra \cdots\ra X_n$$
is an Auslander-Reiten $n$-angle in $\mathcal{S}$ and $X_1,
X_n\notin\add(\oplus_{i=2}^{n-1}X_i)$. Then the two rings
$\End_{\mathcal{S}}(\oplus^{n-1}_{i=1}X_i)/I$ and
$\End_{\mathcal{S}}(\oplus^{n}_{i=2}X_i)/J$ are derived equivalent,
where $I, J$ are defined as in Theorem \ref{Theo}.
\end{Koro}

{\bf Proof.} By \cite[Proposition 3.9]{IY} and Corollary
\ref{coro2}, we can get the conclusion. $\square$

\bigskip

\section{Examples}

In this part, we give an example to illustrate the main result of
this paper.

Consider the 2-representation finite algebra $A$ of type `A'. The
quiver with relation of $A$ is given by the following diagram.
$$
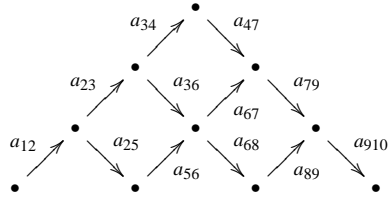
\xy
(0,0)*+{\scriptstyle{\bullet}}="a1",
(8,8)*+{\scriptstyle{\bullet}}="a2",
(16,16)*+{\scriptstyle{\bullet}}="a3",
(24,24)*+{\scriptstyle{\bullet}}="a4",
(16,0)*+{\scriptstyle{\bullet}}="a5",
(24,8)*+{\scriptstyle{\bullet}}="a6",
(32,16)*+{\scriptstyle{\bullet}}="a7",
(32,0)*+{\scriptstyle{\bullet}}="a8",
(40,8)*+{\scriptstyle{\bullet}}="a9",
(48,0)*+{\scriptstyle{\bullet}}="a10", {\ar^{a_{12}} "a1"; "a2"},
{\ar^{a_{23}} "a2"; "a3"}, {\ar^{a_{34}} "a3"; "a4"}, {\ar_{a_{56}}
"a5"; "a6"}, {\ar_{a_{67}} "a6"; "a7"}, {\ar_{a_{89}} "a8"; "a9"},
{\ar^{a_{910}} "a9"; "a10"}, {\ar^{a_{79}} "a7"; "a9"},
{\ar^{a_{47}} "a4"; "a7"}, {\ar^{a_{68}} "a6"; "a8"}, {\ar^{a_{36}}
"a3"; "a6"}, {\ar^{a_{25}} "a2"; "a5"},
\endxy$$
with relations $\{a_{23}a_{36}-a_{25}a_{56},
a_{34}a_{47}-a_{36}a_{67}, a_{67}a_{79}-a_{68}a_{89}, a_{12}a_{25},
a_{56}a_{68}, a_{89}a_{910}\}.$

Assume that $ \nu:=DA\otimes^L_A-: D(A)\ra D(A)$ is the derived
functor of Nakayama functor and $\nu_n=\nu[-n]$. By \cite[Theorem
1]{GKO}, The $2$-cluster tilting subcategory $\mathcal
{U}=\add\{\nu^i_2A\mid i\in \mathbb{Z}\}$ of $D(A)$ is a
$4$-angulated category with suspension functor $\Sigma_4$. And the
Auslander-Reiten quiver of $\mathcal {U}$ is given as follows. (see
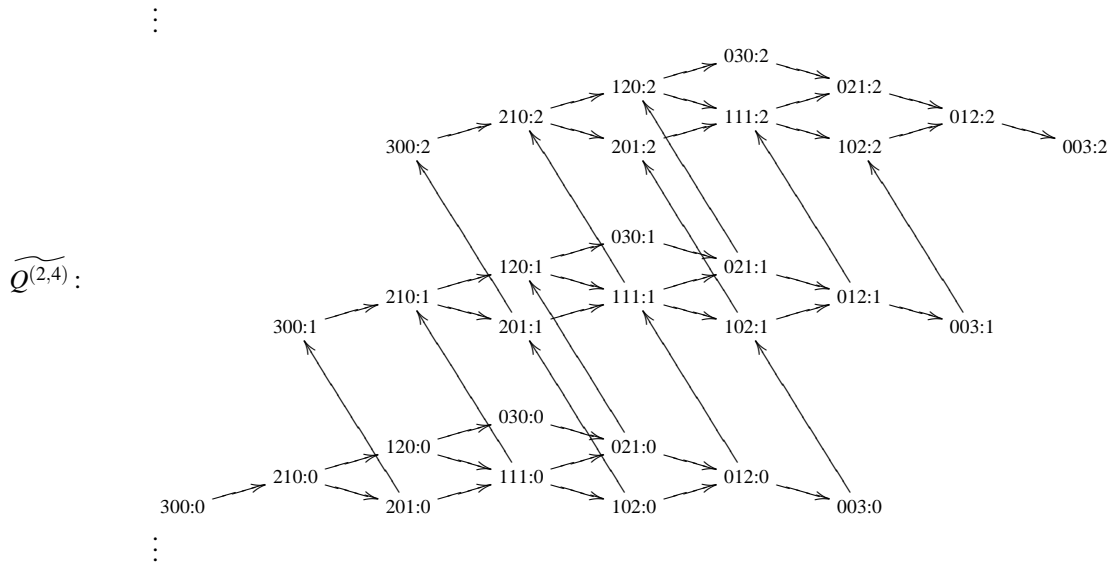
\cite{I2, IO1})
$$\widetilde{Q^{(2,4)}}:\quad\quad\begin{array}{l}
\vdots\\
 \xy (0,0)*+{\scriptstyle{300:0}}="a1",
(15,4)*+{\scriptstyle{210:0}}="a2",
(30,8)*+{\scriptstyle{120:0}}="a3",
(45,12)*+{\scriptstyle{030:0}}="a4",
(30,0)*+{\scriptstyle{201:0}}="a5",
(45,4)*+{\scriptstyle{111:0}}="a6",
(60,8)*+{\scriptstyle{021:0}}="a7",
(60,0)*+{\scriptstyle{102:0}}="a8",
(75,4)*+{\scriptstyle{012:0}}="a9",
(90,0)*+{\scriptstyle{003:0}}="a10", {\ar "a1"; "a2"}, {\ar "a2";
"a3"}, {\ar "a3"; "a4"}, {\ar "a5"; "a6"}, {\ar "a6"; "a7"}, {\ar
"a8"; "a9"}, {\ar "a9"; "a10"}, {\ar "a7"; "a9"}, {\ar "a4"; "a7"},
{\ar "a6"; "a8"}, {\ar "a3"; "a6"}, {\ar "a2"; "a5"},
(15,24)*+{\scriptstyle{300:1}}="b1",
(30,28)*+{\scriptstyle{210:1}}="b2",
(45,32)*+{\scriptstyle{120:1}}="b3",
(60,36)*+{\scriptstyle{030:1}}="b4",
(45,24)*+{\scriptstyle{201:1}}="b5",
(60,28)*+{\scriptstyle{111:1}}="b6",
(75,32)*+{\scriptstyle{021:1}}="b7",
(75,24)*+{\scriptstyle{102:1}}="b8",
(90,28)*+{\scriptstyle{012:1}}="b9",
(105,24)*+{\scriptstyle{003:1}}="b10", {\ar "b1"; "b2"}, {\ar "b2";
"b3"}, {\ar "b3"; "b4"}, {\ar "b5"; "b6"}, {\ar "b6"; "b7"}, {\ar
"b8"; "b9"}, {\ar "b9"; "b10"}, {\ar "b7"; "b9"}, {\ar "b4"; "b7"},
{\ar "b6"; "b8"}, {\ar "b3"; "b6"}, {\ar "b2"; "b5"}, {\ar "a5";
"b1"}, {\ar "a6"; "b2"}, {\ar "a7"; "b3"}, {\ar "a8"; "b5"}, {\ar
"a9"; "b6"}, {\ar "a10"; "b8"}, (30,48)*+{\scriptstyle{300:2}}="c1",
(45,52)*+{\scriptstyle{210:2}}="c2",
(60,56)*+{\scriptstyle{120:2}}="c3",
(75,60)*+{\scriptstyle{030:2}}="c4",
(60,48)*+{\scriptstyle{201:2}}="c5",
(75,52)*+{\scriptstyle{111:2}}="c6",
(90,56)*+{\scriptstyle{021:2}}="c7",
(90,48)*+{\scriptstyle{102:2}}="c8",
(105,52)*+{\scriptstyle{012:2}}="c9",
(120,48)*+{\scriptstyle{003:2}}="c10", {\ar "c1"; "c2"}, {\ar "c2";
"c3"}, {\ar "c3"; "c4"}, {\ar "c5"; "c6"}, {\ar "c6"; "c7"}, {\ar
"c8"; "c9"}, {\ar "c9"; "c10"}, {\ar "c7"; "c9"}, {\ar "c4"; "c7"},
{\ar "c6"; "c8"}, {\ar "c3"; "c6"}, {\ar "c2"; "c5"}, {\ar "b5";
"c1"}, {\ar "b6"; "c2"}, {\ar "b7"; "c3"}, {\ar "b8"; "c5"}, {\ar
"b9"; "c6"}, {\ar "b10"; "c8"},
\endxy\\
\vdots \end{array}$$

Note that the functor $\nu_2$ can be viewed as the automorphism of
$\widetilde{Q^{(2,4)}}$ which send $(l_1,l_2,l_3: i)$ to
$(l_1,l_2,l_3: i-1)$. Select a source map $f_1: 111:0\ra 210:1\oplus
021:0\oplus 102:0$. There is a $4$-angle
$$111:0\stackrel{f_1}\ra 210:1\oplus 021:0\oplus 102:0\ra X_3\stackrel{g}\ra X_4\ra \Sigma_4111:0\quad\quad \quad(*)$$
in $\mathcal {U}$. By \cite[Proposition 3.9]{IY}, $(*)$ is an
Auslander-Reiten $4$-angle in $\mathcal {U}$ and $g$ is a sink map.
By \cite[Theorem 3.10]{IY}, we have $111:0=\nu_2X_4$. Thus,
$$111:0\stackrel{f_1}\ra 210:1\oplus 021:0\oplus 102:0\stackrel{f_2}\ra 120:1\oplus 201:1\oplus 012:0\stackrel{f_3}\ra 111:1\stackrel{f_4}\ra \Sigma_4 111:0$$
is an Auslander-Reiten $4$-angle in $\mathcal {U}$.

We denote $210:1\oplus 021:0\oplus 120:1\oplus102:0\oplus
201:1\oplus012:0$ by $M$. Clearly, the morphism $f_1: 111:0\ra
210:1\oplus 021:0\oplus102:0$ is a left $\add (M)$-approximation of
$111:0$ and the morphism $f_3: 120:1\oplus 201:1\oplus 012:0\ra
111:1$ is a right $\add(M)$-approximation of $111:1$. By Corollary
\ref{nAa}, we can get that the two rings
$\End_{D(\modcat{A})}(111:0\oplus M)/I$ and
$\End_{D(\modcat{A})}(M\oplus 111:1)/J$ are derived equivalent where
$I,J$ are defined as in Theorem \ref{Theo}.

\medskip

\bigskip

{\bf Acknowledgments.} This is part of my Ph.D. thesis, written
under the supervision of Professor Changchang Xi at Beijing Normal
University. I want to express my gratitude to Professor Changchang
Xi for encouragements and useful suggestions.

\medskip \footnotesize{
}

\end{document}